\documentclass[11pt]{article}%
\usepackage{amsmath}
\usepackage{amsthm}
\usepackage{dsfont}
\usepackage{graphicx}
\usepackage{xcolor}
\usepackage{multirow}
\usepackage{booktabs}%
\usepackage{amsfonts}%
\usepackage{amssymb}
\newcommand{\eqnb}{\begin{equation}}
\newcommand{\eqne}{\end{equation}}

\newtheorem{The}{Theorem}

\newtheorem{Cor}[The]{Corollary}
\newtheorem{Lem}{Lemma}

\newtheorem{Rem}{Remark}

\begin{document}

\title{\textbf{Double-End Queues with Non-Poisson Inputs and Their Effective Algorithms}}
\author{Heng-Li Liu$^{a}$, Quan-Lin Li$^{b}$ \thanks{Corresponding author: Q.L. Li
(liquanlin@tsinghua.edu.cn)}, Yan-Xia Chang$^{b}$, Chi Zhang$^{b}$,\\$^{a}$School of Economics and Management Sciences\\Yanshan University, Qinhuangdao 066004, China\\$^{b}$School of Economics and Management\\Beijing University of Technology, Beijing 100124, China}
\maketitle

\begin{abstract}
It is interesting and challenging to study double-ended queues with
First-Come-First-Match discipline under customers' impatient behavior and
non-Poisson inputs. The system stability can be guaranteed by the customers'
impatient behavior, while the existence of impatient customers makes analysis
of such double-ended queues more difficult or even impossible to find an
explicitly analytic solution, thus it becomes more and more important to
develop effective numerical methods in a variety of practical matching
problems. This paper studies a block-structured double-ended queue, whose
block structure comes from two independent Markovian arrival processes (MAPs),
which are non-Poisson inputs. We show that such a queue can be expressed as a
new bilateral quasi birth-and-death (QBD) process which has its own interest.
Based on this, we provide a detailed analysis for both the bilateral QBD
process and the double-ended queue, including the system stability, the queue
size distributions, the average stationary queue lengths, and the sojourn time
of any arriving customers. Furthermore, we develop three effective algorithms
for computing the performance measures (i.e., the probabilities of stationary
queue lengths, the average stationary queue lengths, and the average sojourn
times) of the double-ended queue with non-Poisson inputs. Finally, we use some
numerical examples in tabular and graphical to illustrate how the performance
measures are influenced by some key system parameters. We believe that the
methodology and results described in this paper can be applicable to deal with
more general double-ended queues in practice, and develop some effective
algorithms for the purpose of many actual uses.

\textbf{Keywords:} Double-ended queue; First-Come-First-Match; impatient
customer; Markovian arrival process (MAP); QBD process; RG-factorization.

\end{abstract}

\section{Introduction}

In recent years, we are facing more and more matching problems in practice,
for example, sharing economy, platform service, multilateral market, organ
transplantation, communication network, intelligent manufacturing,
transportation networks, and so on. Based on this, it is interesting but
difficult and challenging to study double-ended queues with
First-Come-First-Match discipline (also known as \textit{matched queues})
under customers' impatient behavior and non-Poisson inputs. Note that the
stability of this system is guaranteed by customers' impatient behavior, while
the existence of impatient customers always makes analysis of such
double-ended queues more difficult or even impossible to find an explicitly
analytic solution, thus it becomes more and more important to develop
effective numerical methods in a variety of practical matching problems. This
motivates us in this paper to develop effective algorithms in a more general
double-ended queue with non-Poisson inputs. Therefore, this paper studies a
block-structured double-ended queue, whose block structure comes from two
independent Markovian arrival processes (MAPs). We show that such a queue can
be expressed as a new bilateral QBD process, and provide a detailed analysis
for the bilateral QBD process, including the system stability, the queue size
distributions, the average stationary queue lengths, and the sojourn time of
any arriving customers. Furthermore, we develop three effective algorithms for
computing the performance measures (i.e., the probabilities of stationary
queue lengths, the average stationary queue lengths, and the average sojourn
times) of the double-ended queue with non-Poisson inputs. Finally, numerical
examples are employed to illustrate how the performance measures are
influenced by key system parameters. We believe that the methodology and
results described in this paper can be applicable to deal with more general
double-ended queues. Also, this can lead to some new theory of bilateral
block-structured Markov processes (e.g., bilateral Markov processes of GI/M/1
type and of M/G/1 type).

Since the double-ended (or matched) queue was first proposed by Kendall
(1951), it has received high attention from many practical applications.
Important examples include \textit{organ transplantation} by Zenios (1999),
Boxma et al. (2011), Stanford et al. (2014) and Elalouf et al. (2018);
\textit{taxi issues} by Giveen (1961, 1963), Kashyap (1965, 1966, 1967), Bhat
(1970), Shi and Lian (2016) and Zhang et al. (2019); \textit{sharing economy}
by Cheng (2016), Sutherland and Jarrahi (2018), Benjaafar and Hu (2019) and
Liu et al. (2021); \textit{transportation} by Browne et al. (1970) and Baik et
al. (2002); \textit{assembly systems} by Hopp and Simon (1989), Som et al.
(1994) and Ramachandran and Delen (2005); \textit{inventory management} by
Sasieni (1961), Porteus (1990), and Axs\"{a}er (2015); \textit{health care} by
Pandey and Gangeshwer (2018); \textit{multimedia synchronization} by Steinmetz
(1990) and Parthasarathy et al. (1999), and so on.

Recently, an emerging hot research topic focuses on ridesharing platform,
which is used to match customers with servers in a bilateral market, for
example, transportation, housing, eating, dressing, and so forth. Plenty of
research has been conducted in this area, including Azevedo and Weyl (2016),
Duenyas et al. (1997), Hu and Zhou (2015), Banerjee and Johari (2019),
Braverman et al. (2019), Liu et al. (2021) and so forth. Meanwhile, many
ridesharing companies in these areas spring up as a result of rapid
development of mobile networks, smart phones and location technologies, such
as Uber in transportation, Airbnb in housing, Eatwith in eating, and Rent the
Runway in dressing. Observing the ridesharing platform reveals that the match
process works as follows: If a customer is in demand of service and makes a
request by using his/her smart phone, then the ridesharing platform will try
to match the customer with a server. Once such a match is successful, the
functioning of the ridesharing platform for this customer is finished
immediately. Therefore, the double-ended queue is an effective mathematical
method for studying ridesharing platforms, as it is a bilateral matching
system or market. Normally, the matching process follows the
First-Come--First-Match principle.

Although the double-ended queue seems simple as it contains only a few random
factors, its analysis is actually very difficult and challenging due to the
fact that the Markov process corresponding to the double-ended queue has a
bidirectional state space $\left\{  \ldots,-2,-1,0,1,2,\ldots\right\}  $.
Therefore, up until now, there is still lack of effective methods for
conducting performance analysis of double-ended queues. In the early study of
double-ended queues, by applying the Markov processes, Sasieni (1961), Giveen
(1961) and Dobbie (1961) established the Chapman-Kolmogorov forward
differential-difference equations with the bidirectional state space. Also,
they introduced customers' impatient behavior to guarantee system stability.
It is well-known that the customers' impatient behavior further makes analysis
of the double-ended queues more difficult and challenging in performance
evaluation of the systems due to the level-dependent structure of the
corresponding Markov processes, which was discussed in Artalejo and
G\'{o}mez-Corral (2008) and Li (2010).

When the two waiting rooms of the double-ended queue are both finite, Jain
(1962) and Kashyap (1965, 1966, 1967) applied the supplementary variable
method to deal with the double-ended queue with a Poisson arrival process and
a renewal arrival process. Takahashi et al. (2000), and Takahashi and
Takahashi (2000) considered a double-ended queue with a Poisson arrival
process and a PH-renewal arrival process. Sharma and Nair (1991) used the
matrix theory to analyze the transient behavior of a double-ended Markovian queue.

When the two waiting rooms are both infinite, Latouche (1981) applied the
matrix-geometric solution to analyze several different bilateral matching
queues with paired input. Conolly et al. (2002) applied the Laplace transform
to discuss the time-dependent performance measures of the double-ended queue
with state-dependent impatience. Di Crescenzo et al. (2012, 2018) discussed
the transient and stationary probability laws of a time-nonhomogeneous
double-ended queue with catastrophes and repairs. Diamant and Baron (2019)
analyzed a double-ended queue with priority and impatient customers, and
derived exact formulae for the stationary queue length distribution and
several useful performance measures. Following the matrix-analytic method
based on the RG-factorizations, Liu et al. (2020) discussed a double-ended
queue with matching batch pair $(m,n)$, which includes two useful examples:
Many-to-one and one-to-many matching requests.

When one waiting room is finite while the other one is infinite, Xu et al.
(1990) discussed a double-ended queue with two Poisson inputs, a PH service
time distribution and a matching proportion $1:r$. They applied the
matrix-geometric solution (see Neuts (1981)) to obtain the stable condition of
the system, and to analyze the stationary queue lengths of both classes of
customers. Further research includes Xu et al. (1993) and Xu and He (1993,
1994). Yuan (1992) applied Markov chains of M/G/1 type (see Neuts (1989)) to
consider a double-ended queue with two Poisson inputs, a general service time
distribution and a matching proportion $1:r$. Li and Cao (1996) further
applied Markov chains of M/G/1 type to deal with a double-ended queue with two
batch Markovian arrival processes (BMAPs), a general service time distribution
and a matching proportion $1:r$. Wu and He (2020) applied the theory of
multi-layer Markov modulated fluid flow (MMFF) processes to analyze a
double-sided queueing model with marked Markovian arrival processes and finite
discrete abandonment times.

In a double-ended queue, if the two classes of customer arrivals are both
general renewal processes, then the ordinary Markov methods (for example,
continuous-time Markov chains, the supplementary variable method, the
matrix-analytic method and so on) will not work well any more. In this case,
the fluid and diffusion approximations become an effective (but approximative)
mathematical method to deal with the more general double-ended queues. Jain
(2000) applied diffusion approximation to discuss the G$^{\text{X}}%
$/G$^{\text{Y}}$/1 double-ended queue. Di Crescenzo et al. (2012, 2018)
discussed a double-ended queue by means of a jump-diffusion approximation. Liu
et al. (2014) discussed diffusion models for the double-ended queues with two
renewal arrival processes. B\"{u}ke and Chen (2017) applied fluid and
diffusion approximations to study the probabilistic matching systems. Liu
(2019) analyzed diffusion approximations for the double-ended queues with
reneging in heavy traffic.

Kim et al. (2010) provided a simulation model for a more general double-ended
queue. Jain (1995) proposed a sample path analysis for the double-ended queue
with time-dependent rates. Af\`{e}che et al. (2014) applied the level-crossing
method to analyze the double-ended batch queue with abandonment.

Hlynka and Sheahan (1987) analyzed the control rates in a double-ended queue
with two Poisson inputs. Gurvich and Ward (2014) discussed dynamic control of
the double-ended queues. B\"{u}ke and Chen (2015) analyzed stabilizing
admission control policies for the probabilistic matching systems. Lee et al.
(2019) discussed optimal control of a time-varying double-ended production
queueing model.

The MAP (Markov arrival process) is a useful mathematical tool, for example,
for describing bursty traffic and dependent arrivals in many practical
systems, such as computer and communication networks, manufacturing systems,
transportation networks and so on. Also, the MAP contains the Poisson process,
the PH-renewal process, and the Markovian Modulated Poisson Process as its
special cases, e.g., see Section 1.5 of Li (2010). Readers may refer to recent
publications for details, among which are Neuts (1979), Chapter 5 in Neuts
(1989), Lucantoni (1991), Narayana and Neuts (1992), Chakravarthy (2001),
Chapter 1 in Li (2010), Cordeiro and Kharoufeh (2010) and references therein.
In the current matching problems (e.g., online bilateral markets, and ride
sharing platform), the arrival data flow is always bursty traffic and
dependent arrivals. Note that the MAP may be a non-renewal process, and it can
express many key features of a practical data flow by means of the statistical
adjustment of multiple parameters. Therefore, it is necessary and useful to
use the MAP inputs to study the double-ended queues.

To be able to deal with level-dependent (or general) Markov processes, the
RG-factorizations were systematically developed in Li (2010). Readers may also
refer to, such as, Li and Cao (2004), Li and Liu (2004) and Li and Zhao (2004)
for early research. Note that the RG-factorizations were successfully applied
to analysis of retrial queues, processor-sharing queues, queues with negative
customers, and queues with impatient customers due to their level-dependent
Markov processes, thus the RG-factorizations play a key role in the study of
double-ended queues whose impatient customers of guaranteeing system stability
directly lead to the level-dependent Markov processes. By using the
RG-factorizations, this paper finds a feasible solution of the double-ended
queue with two MAP inputs and impatient customers, which is more general than
those works in the existing literature. Also, we develop some effective
RG-factorization algorithms which are able to numerically analyze performance
measures of the block-structured double-ended queue.

We summarize the main contributions of this paper as follows:

\begin{itemize}
\item[(1)] We consider a more general block-structured double-ended queue with
non-Poisson (MAP) inputs and impatient customers, and show that such a queue
can be expressed as a new bilateral QBD process which has its own interest. By
using the bilateral QBD process, we can analyze the system stability, the
stationary queue lengths and the sojourn times of this double-ended queue.
Therefore, we develop a new effective RG-factorization method in the study of
double-ended queues, which is different from those in the literature.

\item[(2)] We provide three effective algorithms for computing performance
measures of the block-structured matched queue, in which the bilateral
level-dependent QBD process is decomposed into two unilateral level-dependent
QBD processes that are specifically linked in Level $0$. From the two
unilateral level-dependent QBD processes, we can effectively compute the $R$-,
$U$- and $G$-measures by means of the those approximate algorithms given in
Bright and Taylor (1995, 1997). Based on this, we can numerically compute the
stationary queue lengths. Also, we provide an effective method to discuss the
sojourn times of the block-structured double-ended queue by using the
technique of the first passage times and the PH distributions.

\item[(3)] We use some numerical examples to indicate how the performance
measures of the double-ended queue are influenced by key system parameters. In
addition, the numerical results are also given a simple and interesting
discussion by means of the coupling method of Markov processes.
\end{itemize}

\vskip                                         0.4cm

The structure of this paper is organized as follows. Section 2 describes a
double-ended queue with two MAP inputs and customers' impatient behavior.
Section 3 shows that the double-ended queue can be expressed as a new
bilateral QBD process. By using the bilateral QBD process, we obtain some
stable conditions for the double-ended queues. Section 4 studies the
stationary probability vector of the bilateral QBD process, and compute the
probabilities of stationary queue lengths and the average stationary queue
lengths. Section 5 provides an effective method to discuss the sojourn time of
any arriving customer and to compute the average sojourn time by using the
technique of the first passage times and the PH distributions. Section 6 uses
some numerical examples to indicate how the performance measures of the
double-ended queue are influenced by key system parameters, where three
effective algorithms are developed. Finally, Section 7 gives some concluding remarks.

\section{Model Description}

In this section, we describe a more general block-structured double-ended
queue with two MAP inputs and customers' impatient behavior, and also
introduce operational mechanism, system parameters and basic notation

In the proposed queue, there are two types of customers, called A-customers
and B-customers. Once an A-customer and a B-customer enter their corresponding
buffers, they match each other to constitute a pair and leave the queueing
system immediately, that is, their matching time is zero. Figure 1 provides a
physical illustration for such a double-ended queue.

\begin{figure}[tbh]
\centering         \includegraphics[height=4cm
,width=14cm]{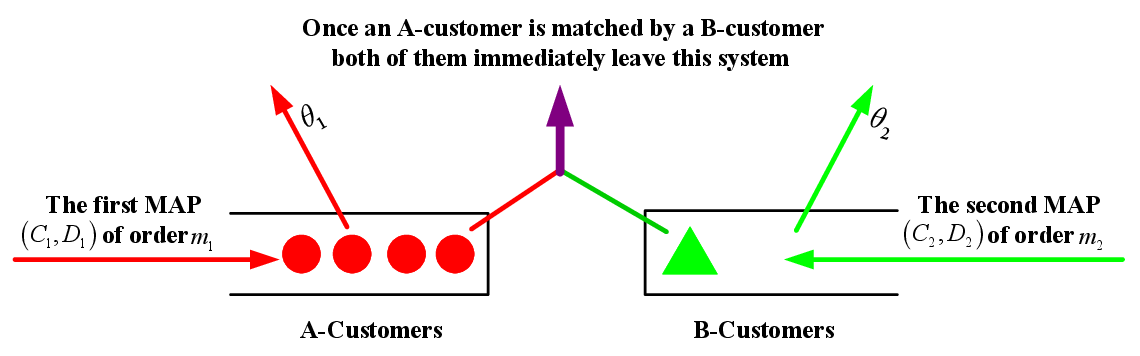}
\caption{A physical illustration of the double-ended queue}%
\end{figure}

Now, we provide a more detailed description for the double-ended queue as follows:

\begin{itemize}
\item[(1)] The A-customers arrive at the queueing system according to a MAP
with irreducible matrix representation $\left(  C_{1},D_{1}\right)  $ of order
$m_{1}$, where $D_{1}\gvertneqq0$, each diagonal element of $C_{1}$ is
negative while its nondiagonal elements are nonnegative, $C_{1}\mathbf{e}%
\lneqq0$, $\left(  C_{1}+D_{1}\right)  \mathbf{e}=0$, and $\mathbf{e}$ is a
column vector of ones with a suitable size. We assume that the Markov process
$C_{1}+D_{1}$ is irreducible and positive recurrent. Let $\alpha_{1}$ be the
stationary probability vector of the Markov process $C_{1}+D_{1}$. Then
$\lambda_{1}=\alpha_{1}D_{1}\mathbf{e}$ is the stationary arrival rate of the
MAP with irreducible matrix representation $\left(  C_{1},D_{1}\right)  $.

Similarly, the B-customers arrive at the queueing system according to a MAP
with irreducible matrix representation $\left(  C_{2},D_{2}\right)  $ of order
$m_{2}$, where $D_{2}\gvertneqq0$, each diagonal element of $C_{2}$ is
negative while its nondiagonal elements are nonnegative, $C_{2}\mathbf{e}%
\lneqq0$, and $\left(  C_{2}+D_{2}\right)  \mathbf{e}=0$. We assume that the
Markov process $C_{2}+D_{2}$ is irreducible and positive recurrent. Let
$\alpha_{2}$ be the stationary probability vector of the Markov process
$C_{2}+D_{2}$. Then $\lambda_{2}=\alpha_{2}D_{2}\mathbf{e}$ is the stationary
arrival rate of the MAP with irreducible matrix representation $\left(
C_{2},D_{2}\right)  $.
\end{itemize}

For the MAPs, readers may refer to Section 1.5 in Chapter 1 of Li
\cite{Li:2010} for a more detailed introduction.

\begin{itemize}
\item[(2)] If an A-customer (resp. a B-customer) stays in the queueing system
for a long time, then she will show some impatience. To capture this
phenomenon, we assume that the impatient time, which is defined as the longest
time that a customer stays in the buffer before leaving, of an A-customer
(resp. a B-customer) is exponentially distributed with impatient rate
$\theta_{1}$ (resp. $\theta_{2}$) for $\theta_{1},\theta_{2}>0$.

\item[(3)] Once an A-customer and a B-customer match as a pair, both of them
immediately leave the queueing system. The matching process follows a
First-Come-First-Match discipline and has the zero matching time. We assume
the waiting spaces of A- and B-customers are all infinite.

\item[(4)] We assume that all the random variables defined above are
independent of each other.
\end{itemize}

\vskip                                            0.4cm

Finally, we provide some necessary interpretation for some practical factors
in the double-ended queue as follows:

(a) The matching times between the A- and B-customers are very short under the
current network environment of data exchange and transmission. Thus the
matching times are regarded as zero so that the matching of A- and B-customers
is complete instantly (or immediately).

(b) The impatient behavior given in Assumption (2) is used to ensure the
stability of the double-ended queue, which can be widely found in the
real-world situations, such as transportation, housing, eating and so forth.

(c) The matching discipline and the zero matching time given in Assumption (3)
indicates that the A- and B-customers cannot simultaneously exist in their
corresponding waiting spaces, while such a phenomenon can also be widely found
in the real-world situations.

\section{A Bilateral QBD Process}

In this section, we show that the block-structured double-ended queue can be
expressed as a new bilateral QBD process with bidirectional infinite sizes. By
using the bilateral QBD process, we obtain some stability conditions of the
double-ended queue.

We denote by $N_{1}\left(  t\right)  $ and $N_{2}\left(  t\right)  $ the
numbers of A- and B-customers in the double-ended queue at time $t\geq0$,
respectively. Let $J_{1}\left(  t\right)  $ and $J_{2}\left(  t\right)  $ be
the phases of two MAPs for the A- and B-customers at time $t$, respectively.
Then the double-ended queue can be modeled as a four-dimensional Markov
process $\left\{  \left(  N_{1}\left(  t\right)  ,J_{1}\left(  t\right)
;N_{2}\left(  t\right)  ,J_{2}\left(  t\right)  \right)  ,t\geq0\right\}  $.

For $t\geq0$, we write all the possible values of $\left(  N_{1}\left(
t\right)  ,N_{2}\left(  t\right)  \right)  $ as
\[
\left\{  \ldots,\left(  0,3\right)  ,\left(  0,2\right)  ,\left(  0,1\right)
,\left(  0,0\right)  ,\left(  1,0\right)  ,\left(  2,0\right)  ,\left(
3,0\right)  ,\ldots\right\}  .
\]
Based on the First-Come-First-Match discipline, it is easy to see that at
least one of the two numbers $N_{1}\left(  t\right)  $ and $N_{2}\left(
t\right)  $ is zero at any time $t\geq0$. Let $N\left(  t\right)
=N_{1}\left(  t\right)  -N_{2}\left(  t\right)  $. Then all the possible
values of $N\left(  t\right)  $ for $t\geq0$ are described as%
\[
\left\{  \ldots,-3,-2,-1,0,1,2,3,\ldots\right\}  .
\]
Thus, the state transition relations of the Markov process $\left\{  \left(
N\left(  t\right)  ,J_{2}\left(  t\right)  ,J_{1}\left(  t\right)  \right)
,\text{ }t\geq0\right\}  $ is depicted in Figure 2. It is easy to see from
Figure 2 that the Markov process $\{(N\left(  t\right)  ,J_{2}\left(
t\right)  $, $J_{1}\left(  t\right)  ),$ $t\geq0\}$ is a new bilateral QBD
process on a state space, given by%
\begin{align*}
\Omega &  =\left\{  \left(  n,i,j\right)  :n\in\left\{  \ldots
,-2,-1,0,1,2,\ldots\right\}  ,\right. \\
&  \left.  i\in\left\{  1,2,\ldots,m_{2}\right\}  ,j\in\left\{  1,2,\ldots
,m_{1}\right\}  \right\}  .
\end{align*}

\begin{figure}[tbh]
\centering
\includegraphics[height=5cm
,width=15cm]{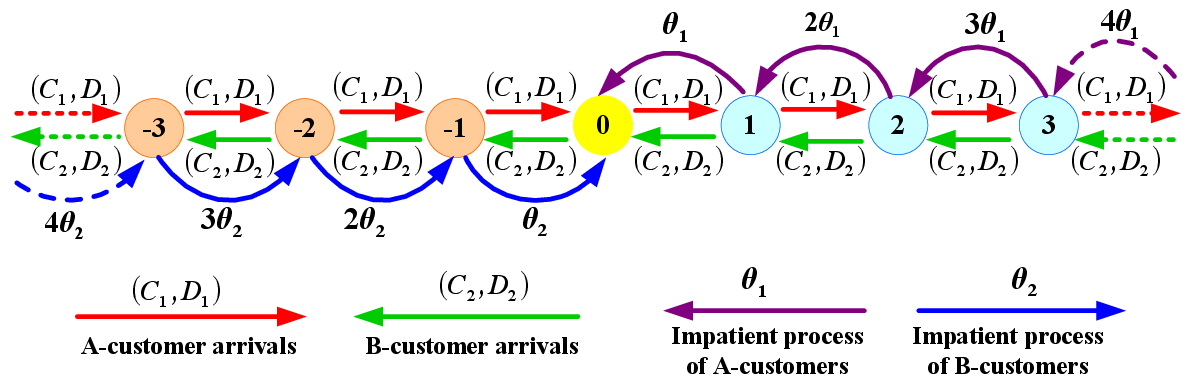} \caption{The state
transition relations of the bilateral QBD process}%
\end{figure}

Based on Figure 2, the infinitesimal generator of the bilateral QBD process
$\{(N\left(  t\right)  ,J_{2}\left(  t\right)  $, $J_{1}\left(  t\right)  ),$
$t\geq0\}$ is given by
\begin{equation}
Q=\left(
\begin{array}
[c]{ccccccccccc}%
\ddots & \ddots & \ddots &  &  &  &  &  &  &  & \\
& B_{0}^{\left(  -3\right)  } & B_{1}^{\left(  -3\right)  } & B_{2}^{\left(
-3\right)  } &  &  &  &  &  &  & \\
&  & B_{0}^{\left(  -2\right)  } & B_{1}^{\left(  -2\right)  } &
B_{2}^{\left(  -2\right)  } &  &  &  &  &  & \\
&  &  & B_{0}^{\left(  -1\right)  } & B_{1}^{\left(  -1\right)  } &
\fbox{$B_{2}^{\left(  -1\right)  }$} &  &  &  &  & \\
&  &  &  & \fbox{$B_{0}^{\left(  0\right)  }$} & \fbox{$B_{1}^{\left(
0\right)  }+A_{1}^{\left(  0\right)  }$} & \fbox{$A_{0}^{\left(  0\right)  }$}
&  &  &  & \\
&  &  &  &  & \fbox{$A_{2}^{\left(  1\right)  }$} & A_{1}^{\left(  1\right)  }
& A_{0}^{\left(  1\right)  } &  &  & \\
&  &  &  &  &  & A_{2}^{\left(  2\right)  } & A_{1}^{\left(  2\right)  } &
A_{0}^{\left(  2\right)  } &  & \\
&  &  &  &  &  &  & A_{2}^{\left(  3\right)  } & A_{1}^{\left(  3\right)  } &
A_{0}^{\left(  3\right)  } & \\
&  &  &  &  &  &  &  & \ddots & \ddots & \ddots
\end{array}
\right)  \label{equ1}%
\end{equation}
where%
\[
A_{1}^{\left(  0\right)  }=I\otimes C_{1},\text{ \ }B_{1}^{\left(  0\right)
}=C_{2}\otimes I,
\]%
\[
A_{0}^{\left(  k\right)  }=I\otimes D_{1},\text{ \ }k\geq0,
\]%
\[
A_{1}^{\left(  k\right)  }=C_{2}\oplus C_{1}-k\theta_{1}I,\text{ }%
A_{2}^{\left(  k\right)  }=D_{2}\otimes I+k\theta_{1}I,\text{ \ }k\geq1,
\]%
\[
B_{0}^{\left(  k\right)  }=D_{2}\otimes I,\text{ \ }k\leq0,
\]%
\[
B_{1}^{\left(  k\right)  }=C_{2}\oplus C_{1}+k\theta_{2}I,\text{ }%
B_{2}^{\left(  k\right)  }=I\otimes D_{1}-k\theta_{2}I,\text{ \ }k\leq-1.
\]
where $\oplus$ and $\otimes$ are the Kronecker sum and the Kronecker product
for two matrices, respectively.

Now, the key is how to discuss the new bilateral QBD process $Q$, which is an
interesting issue in the study of stochastic models. Note that such a
bilateral QBD process was first introduced in Latouche (1981) and further
analyzed in Li and Cao (2004).

To analyze the double-ended queue, it is observed that Level $0$, together
with the suitable block decomposition $B_{1}^{\left(  0\right)  }%
+A_{1}^{\left(  0\right)  }$, plays a key role, where%
\begin{align*}
\text{Level }0=  &  \left\{  \left(  0,1,1\right)  ,\left(  0,1,2\right)
,\ldots,\left(  0,1,m_{1}\right)  ;\left(  0,2,1\right)  ,\left(
0,2,2\right)  ,\ldots,\left(  0,2,m_{1}\right)  ;\right. \\
&  \left.  \ldots;\left(  0,m_{2},1\right)  ,\left(  0,m_{2},2\right)
,\ldots,\left(  0,m_{2},m_{1}\right)  \right\}  .
\end{align*}
By using the block decomposition $B_{1}^{\left(  0\right)  }+A_{1}^{\left(
0\right)  }$, we divide the bilateral QBD process $Q$ into two unilateral QBD
processes: $Q_{A}$ and $Q_{B}$, both of which are interlinked from two
different state space directions (upward and downward, respectively) by means
of Level $0$ with the blocks $A_{1}^{\left(  0\right)  }$ and $A_{0}^{\left(
0\right)  }$, and the blocks $B_{1}^{\left(  0\right)  }$ and $B_{0}^{\left(
0\right)  }$, respectively. Thus the infinitesimal generators of the two
unilateral QBD processes are respectively given by
\[
Q_{A}=\left(
\begin{array}
[c]{ccccccc}%
A_{1}^{\left(  0\right)  } & A_{0}^{\left(  0\right)  } &  &  &  &  & \\
A_{2}^{\left(  1\right)  } & A_{1}^{\left(  1\right)  } & A_{0}^{\left(
1\right)  } &  &  &  & \\
& A_{2}^{\left(  2\right)  } & A_{1}^{\left(  2\right)  } & A_{0}^{\left(
2\right)  } &  &  & \\
&  & A_{2}^{\left(  3\right)  } & A_{1}^{\left(  3\right)  } & A_{0}^{\left(
3\right)  } &  & \\
&  &  & A_{2}^{\left(  4\right)  } & A_{1}^{\left(  4\right)  } &
A_{0}^{\left(  4\right)  } & \\
&  &  &  & \ddots & \ddots & \ddots
\end{array}
\right)
\]
and%
\[
Q_{B}=\left(
\begin{array}
[c]{ccccccc}%
B_{1}^{\left(  0\right)  } & B_{0}^{\left(  0\right)  } &  &  &  &  & \\
B_{2}^{\left(  -1\right)  } & B_{1}^{\left(  -1\right)  } & B_{0}^{\left(
-1\right)  } &  &  &  & \\
& B_{2}^{\left(  -2\right)  } & B_{1}^{\left(  -2\right)  } & B_{0}^{\left(
-2\right)  } &  &  & \\
&  & B_{2}^{\left(  -3\right)  } & B_{1}^{\left(  -3\right)  } &
B_{0}^{\left(  -3\right)  } &  & \\
&  &  & B_{2}^{\left(  -4\right)  } & B_{1}^{\left(  -4\right)  } &
B_{0}^{\left(  -4\right)  } & \\
&  &  &  & \ddots & \ddots & \ddots
\end{array}
\right)  .
\]

\vskip                                               0.4cm

In the remainder of this section, we analyze the stability of the double-ended
queue by means of that of the two unilateral QBD processes $Q_{A}$ and $Q_{B}%
$, because the bilateral QBD process $Q$ can be divided into the two
unilateral QBD processes with $Q_{A}$ and $Q_{B}$.

The following lemma provides useful relations of stability among the bilateral
QBD process $Q$ and the two unilateral QBD processes $Q_{A}$ and $Q_{B}$. Its
proof is quite straightforward and is omitted here for brevity.

\begin{Lem}
To study the stability of the bilateral QBD process $Q$, we have

(a) The bilateral QBD process $Q$ is positive recurrent if the two QBD
processes $Q_{A}$ and $Q_{B}$ are both positive recurrent.

(b) The bilateral QBD process $Q$ is null recurrent if the QBD process $Q_{A}$
is recurrent and the QBD process $Q_{B}$ is null recurrent; or the QBD process
$Q_{A}$ is null recurrent and the QBD process $Q_{B}$ is recurrent.

(c) The bilateral QBD process $Q$ is transient if at least one of the two QBD
processes $Q_{A}$ and $Q_{B}$ is transient.
\end{Lem}

It is worth noting that the impatience of the two classes of customers plays a
key role in guaranteeing the stability of the double-ended queue (or bilateral
QBD process), in which we consider three different cases: $\left(  \theta
_{1},\theta_{2}\right)  >0$; $\left(  \theta_{1},\theta_{2}\right)  =0$; and
either $\theta_{1}>0,\theta_{2}=0$ or $\theta_{1}=0,\theta_{2}>0$. From the
three cases, we use Lemma 1 to conduct some simple analysis.

The following theorem provides a necessary and sufficient condition for the
stability of the block-structured double-ended queue with $\left(  \theta
_{1},\theta_{2}\right)  >0$.

\begin{The}
\ \label{The1}If $\left(  \theta_{1},\theta_{2}\right)  >0$, then the
bilateral QBD process $Q$ must be irreducible and positive recurrent. Thus the
block-structured double-ended queue is stable.
\end{The}

\textbf{Proof.} \textit{ } Please see (a) Proof of Theorem \ref{The1} in the appendix.

\begin{Rem}
Theorem \ref{The1} shows that the stability of the block-structured
double-ended queue depends on only the two impatient rates $\theta_{1}%
,\theta_{2}>0$, and has nothing to do with the two MAP inputs. This is true
through comparing the rate of the upward shift with the rate of the downward
shift from some levels, and further by using the mean drift method for the
stability of Markov processes, e.g., see Li \cite{Li:2010}.
\end{Rem}

In what follows we consider the influence of the two impatient rates
$\theta_{1}$ and $\theta_{2}$ on the stability of the system.

For $\theta_{1}=\theta_{2}=0$, we discuss the stable conditions of the
bilateral QBD process related to the double-ended queue.

\begin{Cor}
\label{Cor:Proof}Suppose $\theta_{1}=\theta_{2}=0$.

(1) If $\lambda_{1}\neq\lambda_{2}$, then the bilateral QBD process $Q$ is
transient, and the double-ended queue is also transient.

(2) If $\lambda_{1}=\lambda_{2}$, then the bilateral QBD process $Q$ is null
recurrent, and the double-ended queue is also null recurrent.
\end{Cor}

\textbf{Proof.}\textit{ }(1) Since $\lambda_{1}\neq\lambda_{2}$, we shown only
the case with $\lambda_{1}>\lambda_{2}$; while the case with $\lambda
_{1}<\lambda_{2}$ can be dealt with similarly.

Note that, when $\theta_{1}=\theta_{2}=0$ and $\lambda_{1}>\lambda_{2}$, it is
easy to see that the QBD process $Q_{A}$ is transient and the QBD process
$Q_{B}$ is positive recurrent. Thus the bilateral QBD process $Q$ is
transient, so that the double-ended queue is also transient.

(2) Note that, when $\theta_{1}=\theta_{2}=0$ and $\lambda_{1}=\lambda_{2}$,
it is easy to see that the QBD processes $Q_{A}$ and $Q_{B}$ are both null
recurrent. Thus the bilateral QBD process $Q$ is null recurrent, so that the
double-ended queue is also null recurrent. This completes the proof.
\hfill$\blacksquare$

Finally, we consider the case with either $\theta_{1}>0$ and $\theta_{2}=0$ or
$\theta_{1}=0$ and $\theta_{2}>0$.

The following corollary only considers the case with $\theta_{1}>0$ and
$\theta_{2}=0$; while another case with $\theta_{1}=0$ and $\theta_{2}>0$ can
be analyzed similarly and is omitted for brevity.

\begin{Cor}
\label{Cor:Imp}Suppose $\theta_{1}>0$ and $\theta_{2}=0$.

(1) If $\lambda_{1}=\lambda_{2}$, then the bilateral QBD process $Q$ is null
recurrent, and the double-ended queue is also null recurrent.

(2) If $\lambda_{1}>\lambda_{2}$, then the bilateral QBD process $Q$ is
positive recurrent, and the double-ended queue is also positive recurrent.

(3) If $\lambda_{1}<\lambda_{2}$, then the bilateral QBD process $Q$ is
transient, and the double-ended queue is also transient.
\end{Cor}

\textbf{Proof.}\textit{ }(1) If $\lambda_{1}=\lambda_{2}$, $\theta_{1}>0$ and
$\theta_{2}=0$, then the QBD process $Q_{A}$ is positive recurrent while the
QBD process $Q_{B}$ is null recurrent. Thus, the bilateral QBD process $Q$ is
null recurrent, so that the double-ended queue is also null recurrent.

(2) If $\lambda_{1}>\lambda_{2}$, $\theta_{1}>0$ and $\theta_{2}=0$, then the
QBD processes $Q_{A}$ and $Q_{B}$ are both positive recurrent. Thus, the
bilateral QBD process $Q$ is positive recurrent, so that the double-ended
queue is also positive recurrent.

(3) If $\lambda_{1}<\lambda_{2}$, $\theta_{1}>0$ and $\theta_{2}=0$, then the
QBD process $Q_{A}$ is positive recurrent while the QBD process $Q_{B}$ is
transient. Thus, the bilateral QBD process $Q$ is transient, so that the
double-ended queue is also transient. This completes the proof. \hfill
$\blacksquare$

\begin{Rem}
From Corollary \ref{Cor:Proof}, it is seen that if $\theta_{1}=\theta_{2}=0$,
then the double-ended queue can not be positive recurrent. To guarantee the
stability of the double-ended queue, we must introduce the impatient
customers, i.e., $\theta_{1}>0$ or $\theta_{2}>0$. This is given a detailed
analysis in Theorem \ref{The1} for $\theta_{1}>0$ and $\theta_{2}>0$; and
Corollary \ref{Cor:Imp} for either $\theta_{1}>0$ and $\theta_{2}=0$, or
$\theta_{1}=0$ and $\theta_{2}>0$.
\end{Rem}

\section{The Stationary Queue Length}

In this section, we first provide a matrix-product expression for the
stationary probability vector of the bilateral QBD process by means of the
RG-factorizations. Then we provide performance analysis of the
block-structured double-ended queue.

We write%
\[
p_{k;i,j}\left(  t\right)  =P\left\{  N\left(  t\right)  =k;\text{ \ }%
J_{2}\left(  t\right)  =i,J_{1}\left(  t\right)  =j\right\}  .
\]
Since the bilateral QBD process is stable, we have
\[
\pi_{k;i,j}=\lim_{t\rightarrow+\infty}p_{k;i,j}\left(  t\right)  .
\]
For any integer $k=\ldots,-2,-1,0,1,2,\ldots$, we write%
\[
\pi_{k}=\left(  \pi_{k;1,1},\pi_{k;1,2},\ldots,\pi_{k;1,m_{1}};\pi_{k;2,1}%
,\pi_{k;2,2},\ldots,\pi_{k;2,m_{1}};\ldots;\pi_{k;m_{2},1},\pi_{k;m_{2}%
,2},\ldots,\pi_{k;m_{2},m_{1}}\right)
\]
and%
\[
\pi=\left(  \ldots,\pi_{-2},\pi_{-1},\pi_{0},\pi_{1},\pi_{2},\ldots\right)  .
\]

To compute the stationary probability vector of the bilateral QBD process $Q$,
we first need to compute the stationary probability vectors of the two
unilateral QBD processes $Q_{A}$ and $Q_{B}$. Then we use Levels $1$, $0$ and
$-1$ to determine the stationary probability vectors on the interaction
boundary (i.e., Levels $1$, $0$ and $-1$) of the bilateral QBD process $Q$.

Note that the two unilateral QBD processes $Q_{A}$ and $Q_{B}$ are
level-dependent, thus we need to apply the RG-factorizations given in Li
(2010) to calculate their stationary probability vectors. To this end, we need
to introduce the UL-type $U$-, $R$\textit{-} and $G$\textit{-}measures for the
two unilateral QBD processes $Q_{A}$ and $Q_{B}$, respectively. In fact, an
early analysis for such a level-dependent QBD process was given in Ramaswami
and Taylor (1996) and Li and Cao (2004).

For the unilateral QBD process $Q_{A}$, we define the UL-type $U$-,
$R$\textit{-} and $G$\textit{-}measures as%
\begin{equation}
U_{k}=A_{1}^{\left(  k\right)  }+A_{0}^{\left(  k\right)  }\left(
-U_{k+1}^{-1}\right)  A_{2}^{\left(  k+1\right)  },\text{ \ }k\geq0,
\label{equ2}%
\end{equation}%
\[
R_{k}=A_{0}^{\left(  k\right)  }\left(  -U_{k+1}^{-1}\right)  ,\text{ \ }%
k\geq0,
\]
and%
\[
G_{k}=\left(  -U_{k}^{-1}\right)  A_{2}^{\left(  k\right)  },\text{ \ }%
k\geq1.
\]

On the other hand, it is well-known from Ramaswami and Taylor (1996) that the
matrix sequence $\left\{  R_{k},k\geq0\right\}  $ is the minimal nonnegative
solution to the system of nonlinear matrix equations%
\begin{equation}
A_{0}^{\left(  k\right)  }+R_{k}A_{1}^{\left(  k+1\right)  }+R_{k}R_{k+1}%
A_{2}^{\left(  k+2\right)  }=0,\text{ \ }k\geq0, \label{equ3}%
\end{equation}
and the matrix sequence $\left\{  G_{k},k\geq1\right\}  $ is the minimal
nonnegative solution to the system of nonlinear matrix equations%
\begin{equation}
A_{0}^{\left(  k\right)  }G_{k+1}G_{k}+A_{1}^{\left(  k\right)  }G_{k}%
+A_{2}^{\left(  k\right)  }=0,\text{ \ }k\geq1. \label{equ4}%
\end{equation}
Once the matrix sequence $\left\{  R_{k},k\geq0\right\}  $ or $\left\{
G_{k},k\geq1\right\}  $ is given, for $k\geq0$ we have%
\begin{align*}
U_{k}  &  =A_{1}^{\left(  k\right)  }+A_{0}^{\left(  k\right)  }\left(
-U_{k+1}^{-1}\right)  A_{2}^{\left(  k+1\right)  }\\
&  =A_{1}^{\left(  k\right)  }+R_{k}A_{2}^{\left(  k+1\right)  }\\
&  =A_{1}^{\left(  k\right)  }+A_{0}^{\left(  k\right)  }G_{k+1}.
\end{align*}

For the unilateral QBD process $Q_{A}$, by following the method described in
Chapter 1 of Li (2010) or Li and Cao (2004), the UL-type RG-factorization is
given by%
\begin{equation}
Q_{A}=\left(  I-R_{U}\right)  U_{D}\left(  I-G_{L}\right)  , \label{equ5}%
\end{equation}
where%
\[
U_{D}=\text{diag}\left(  U_{0},U_{1},U_{2},U_{3},\ldots\right)  ,
\]%
\[
R_{U}=\left(
\begin{array}
[c]{ccccc}%
0 & R_{0} &  &  & \\
& 0 & R_{1} &  & \\
&  & 0 & R_{2} & \\
&  &  & 0 & \ddots\\
&  &  &  & \ddots
\end{array}
\right)  ,\text{ \ \ \ }G_{L}=\left(
\begin{array}
[c]{ccccc}%
0 &  &  &  & \\
G_{1} & 0 &  &  & \\
& G_{2} & 0 &  & \\
&  & G_{3} & 0 & \\
&  &  & \ddots & \ddots
\end{array}
\right)  .
\]

Let $\pi_{A}=\left(  \pi_{0}^{A},\pi_{1}^{A},\pi_{2}^{A},\ldots\right)  $ be
the stationary probability vector of the unilateral QBD process $Q_{A}$. Then
by applying the UL-type RG-factorization and using the $R$-measure $\left\{
R_{k}:k\geq1\right\}  $, we have
\begin{equation}
\pi_{k}^{A}=\pi_{1}^{A}R_{1}R_{2}\cdots R_{k-1},\text{ \ }k\geq2. \label{equ6}%
\end{equation}

By conducting a similar analysis to that shown above, we can set up the
stationary probability vector, $\pi_{B}=\left(  \pi_{0}^{B},\pi_{1}^{B}%
,\pi_{2}^{B},\ldots\right)  $, of the unilateral QBD process $Q_{B}$. Here, we
only provide the $R$-measure $\left\{  \mathbb{R}_{k}:k\leq0\right\}  $, while
the $U$-measure $\left\{  \mathbb{U}_{k}:k\leq0\right\}  $ and $G$-measure
$\left\{  \mathbb{G}_{k}:k\leq-1\right\}  $ can be given easily and is omitted
for brevity.

Let the matrix sequence $\left\{  \mathbb{R}_{k},k\leq0\right\}  $ be the
minimal nonnegative solution to the system of nonlinear matrix equations%
\begin{equation}
B_{0}^{\left(  k\right)  }+\mathbb{R}_{k}B_{1}^{\left(  k-1\right)
}+\mathbb{R}_{k}\mathbb{R}_{k-1}B_{2}^{\left(  k-2\right)  }=0,\text{ }k\leq0,
\label{equ7}%
\end{equation}
By using the $R$-measure $\left\{  \mathbb{R}_{k}:k\leq-1\right\}  $, we
obtain
\begin{equation}
\pi_{k}^{B}=\pi_{-1}^{B}\mathbb{R}_{-1}\mathbb{R}_{-2}\cdots\mathbb{R}%
_{k+1},\text{ \ }k\leq-2. \label{equ8}%
\end{equation}

\begin{Rem}
In the above analysis, one of the main purposes of applying the UL-type
RG-factorization is to show that our computational procedures can be easily
extended and generalized to deal with more general bilateral block-structured
Markov processes, such as bilateral Markov processes of M/G/1 type, bilateral
Markov processes of GI/M/1 type and so on. See Chapter 1 in Li (2010) for more details.
\end{Rem}

The following theorem expresses the stationary probability vector $\pi=\left(
\ldots,\pi_{-2},\pi_{-1},\right.  $

\noindent$\left.  \pi_{0},\pi_{1},\pi_{2},\ldots\right)  $ of the bilateral
QBD process $Q$ in terms of the stationary probability vectors: $\pi
_{A}=\left(  \pi_{2}^{A},\pi_{3}^{A},\pi_{4}^{A},\ldots\right)  $ given in
(\ref{equ6}), and $\pi_{B}=\left(  \pi_{-2}^{B},\pi_{-3}^{B},\pi_{-4}%
^{B},\ldots\right)  $ given in (\ref{equ8}).

\begin{The}
\label{The2}The stationary probability vector $\pi$ of the bilateral QBD
process $Q$ is given by
\begin{equation}
\pi_{k}=c\widetilde{\pi}_{k}, \label{equ10}%
\end{equation}
and%
\begin{equation}
\widetilde{\pi}_{k}=\left\{
\begin{array}
[c]{l}%
\widetilde{\pi}_{-1}\mathbb{R}_{-1}\mathbb{R}_{-2}\cdots\mathbb{R}%
_{k+1},\text{ }k\leq-2,\\
\widetilde{\pi}_{-1},\\
\widetilde{\pi}_{0},\\
\widetilde{\pi}_{1},\\
\widetilde{\pi}_{1}R_{1}R_{2}\cdots R_{k-1},\text{ \ }k\geq2,
\end{array}
\right.  \label{equ11}%
\end{equation}
where the three vectors $\widetilde{\pi}_{-1},\widetilde{\pi}_{0}%
,\widetilde{\pi}_{1}$ are uniquely determined by the following system of
linear equations%
\begin{equation}
\left\{
\begin{array}
[c]{l}%
\widetilde{\pi}_{-1}\mathbb{R}_{-1}B_{2}^{\left(  -2\right)  }+\widetilde{\pi
}_{-1}B_{1}^{\left(  -1\right)  }+\widetilde{\pi}_{0}B_{0}^{\left(  0\right)
}=0,\\
\widetilde{\pi}_{-1}B_{2}^{\left(  -1\right)  }+\widetilde{\pi}_{0}\left(
B_{1}^{\left(  0\right)  }+A_{1}^{\left(  0\right)  }\right)  +\widetilde{\pi
}_{1}A_{2}^{\left(  1\right)  }=0,\\
\widetilde{\pi}_{0}A_{0}^{\left(  0\right)  }+\widetilde{\pi}_{1}%
A_{1}^{\left(  1\right)  }+\widetilde{\pi}_{1}R_{1}A_{2}^{\left(  2\right)
}=0,
\end{array}
\right.  \label{equ12}%
\end{equation}
and the positive constant $c$ is uniquely given by%
\begin{equation}
c=\frac{1}{\sum_{k\leq-2}\widetilde{\pi}_{-1}\mathbb{R}_{-1}\mathbb{R}%
_{-2}\cdots\mathbb{R}_{k+1}\mathbf{e}+\widetilde{\pi}_{-1}\mathbf{e}%
+\widetilde{\pi}_{0}\mathbf{e}+\widetilde{\pi}_{1}\mathbf{e}+\sum
_{k=2}^{\infty}\widetilde{\pi}_{1}R_{1}R_{2}\cdots R_{k-1}\mathbf{e}}.
\label{equ121}%
\end{equation}
\end{The}

\textbf{Proof.}\textit{ }Please see (a) Proof of Theorem \ref{The2} in the appendix.

\begin{Rem}
(a) Just like Markovian retrial queues, analysis of queues with impatient
customers is in the face of level-dependent Markov processes, hence it is
difficult and challenging to have an explicit analytical expression for their
stationary probabilities.

(b) The retrial customers change the arrival process to be state-dependent;
while the impatient customers make the service process to be state-dependent.
Thus both of them lead to the level-dependent Markov processes whose analysis
is completed by means of almost the only way: RG-factorizations. See Chapter 2
of Li (2010) and some remarks in Chapter 5 of Artalejo and G\'{o}mez-Corral
(2008). Also, the stationary probability vector is always expressed as the
matrix-product solution.
\end{Rem}

In the remainder of this section, we provide performance analysis of the
block-structured double-ended queue by means of the stationary probability
vector of the bilateral QBD process.

Note that the block-structured double-ended queue must be stable for $\left(
\theta_{1},\theta_{2}\right)  >0$. Thus, we denote by $\mathcal{Q}%
,\mathcal{Q}^{\left(  1\right)  }$ and $\mathcal{Q}^{\left(  2\right)  }$ the
stationary queue lengths of the double-ended queue, the A- and B-customers, respectively.

By using Theorem \ref{The2}, we can provide some useful performance measures
as follows:

\textbf{(a) A stationary queue length is zero}

The stationary probability that there is no A-customer is given by%
\[
P\left\{  \mathcal{Q}^{\left(  1\right)  }=0\right\}  =\sum_{k\leq0}\pi
_{k}\mathbf{e}.
\]

The stationary probability that there is no B-customer is given by%
\[
P\left\{  \mathcal{Q}^{\left(  2\right)  }=0\right\}  =\sum_{k=0}^{\infty}%
\pi_{k}\mathbf{e}.
\]

The stationary probability that there is neither A-customer nor B-customer is
given by%
\[
P\left\{  \mathcal{Q}=0\right\}  =P\left\{  \mathcal{Q}^{\left(  1\right)
}=0\text{ and }\mathcal{Q}^{\left(  2\right)  }=0\right\}  =\pi_{0}%
\mathbf{e}.
\]

\textbf{(b) The average stationary queue lengths}
\begin{align*}
E\left[  \mathcal{Q}\right]  =  &  E\left[  \mathcal{Q}\text{ }|\text{
}\mathcal{Q}^{\left(  1\right)  }>0,\mathcal{Q}^{\left(  2\right)  }=0\right]
P\left\{  \mathcal{Q}^{\left(  1\right)  }>0\right\} \\
&  +E\left[  \mathcal{Q}\text{ }|\text{ }\mathcal{Q}^{\left(  2\right)
}>0,\mathcal{Q}^{\left(  1\right)  }=0\right]  P\left\{  \mathcal{Q}^{\left(
2\right)  }>0\right\} \\
=  &  E\left[  \mathcal{Q}^{\left(  1\right)  }\right]  P\left\{
\mathcal{Q}^{\left(  1\right)  }>0\right\}  +E\left[  \mathcal{Q}^{\left(
2\right)  }\right]  P\left\{  \mathcal{Q}^{\left(  2\right)  }>0\right\}  ,
\end{align*}
where%
\begin{align*}
E\left[  \mathcal{Q}^{\left(  1\right)  }\right]   &  =\sum_{k=1}^{\infty}%
k\pi_{k}\mathbf{e}=\sum_{k=1}^{\infty}kc\widetilde{\pi}_{k}\mathbf{e}\\
&  =c\widetilde{\pi}_{1}\mathbf{e+}\sum_{k=2}^{\infty}kc\widetilde{\pi}%
_{1}R_{1}R_{2}\cdots R_{k-1}\mathbf{e,}%
\end{align*}%
\[
P\left\{  \mathcal{Q}^{\left(  1\right)  }>0\right\}  =\sum_{k=1}^{\infty}%
\pi_{k}\mathbf{e=}c\widetilde{\pi}_{1}\mathbf{e+}\sum_{k=2}^{\infty
}c\widetilde{\pi}_{1}R_{1}R_{2}\cdots R_{k-1}\mathbf{e,}%
\]
and%
\begin{align*}
E\left[  \mathcal{Q}^{\left(  2\right)  }\right]   &  =\sum_{k\leq-1}\left(
-k\right)  \pi_{k}\mathbf{e}=\sum_{k\leq-1}\left(  -k\right)  c\widetilde{\pi
}_{k}\mathbf{e}\\
&  =c\widetilde{\pi}_{-1}\mathbf{e+}\sum_{k\leq-2}\left(  -k\right)
c\widetilde{\pi}_{-1}\mathbb{R}_{-1}\mathbb{R}_{-2}\cdots\mathbb{R}%
_{k+1}\mathbf{e,}%
\end{align*}%
\[
P\left\{  \mathcal{Q}^{\left(  2\right)  }>0\right\}  =\sum_{k\leq-1}\pi
_{k}\mathbf{e}=c\widetilde{\pi}_{-1}\mathbf{e+}\sum_{k\leq-2}c\widetilde{\pi
}_{-1}\mathbb{R}_{-1}\mathbb{R}_{-2}\cdots\mathbb{R}_{k+1}\mathbf{e.}%
\]

\section{The Sojourn Time}

In this section, we provide an effective method for analyzing the sojourn time
of any arriving customer, and for computing the average sojourn time.

Note that analysis of the sojourn times is symmetrical and similar in the
double-ended queue, thus our discussion mainly focuses on the sojourn time of
an arriving A-customer, while that of an arriving B-customer can be dealt with similarly.

In the double-ended queue, the sojourn time is the time interval from the
arrival epoch of a customer to its departure time. Let $W_{A}^{\left(
k\right)  }$ be the sojourn time of an arriving A-customer that there are
$k-1$ A-customers in front of her at the moment of her arrival for $k\geq1$.
For convenience of computation, we assume that the arrival moment of this
A-customer is time $0$.

If the arriving A-customer finds $k-1$ A-customers in front of her at time $0$
(i.e., the moment of her arrival), then we denote by $N_{1}\left(  t\right)  $
the number of A-customers in the system at time $t\geq0$, and by $M\left(
t\right)  $ the number of A-customers waiting for the matching service in
front of the arriving A-customer at time $t$. Let $J_{2}\left(  t\right)  $ be
the phase of the MAP of the B-customers at time $t$.

From the model descriptions, it is easy to see that $\left\{  \left(
N_{1}\left(  t\right)  ,M\left(  t\right)  ,J_{2}\left(  t\right)  \right)
,t\geq0\right\}  $ is a Markov process with an absorption state $\Delta$,
where $\Delta$ denotes such a state that this arriving A-customer leaves the
double-ended queue. Thus, the state transition relations of the Markov process
are depicted in Figure 3, in which $J_{2}\left(  t\right)  $ are ignored for
the sake of simplicity.

It is obvious that the sojourn time $W_{A}^{\left(  k\right)  }$ is the first
passage time that the Markov process $\left\{  \left(  N_{1}\left(  t\right)
,M\left(  t\right)  ,J_{2}\left(  t\right)  \right)  ,t\geq0\right\}  $
arrives at the absorption state $\Delta$ for the first time. Note that the
Markov process starts at state $\left(  k,k-1\right)  $ at time $0$.

From Figure 3, the state space of the Markov process with the absorption state
$\Delta$ is given by%
\[
\Xi=\left\{  \Delta\right\}  \cup\left\{  \left(  k-i,k-i-1\right)
:i=0,1,\ldots,k-1\right\}  .
\]

\begin{figure}[tbh]
\centering             \includegraphics[height=5cm, width=15cm]{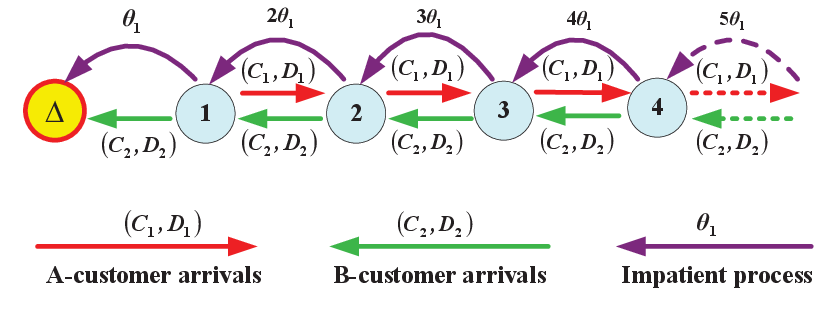}
\caption{The state transition relations of the Markov process with the
absorption state $\Delta$}%
\end{figure}

From Figure 3 and the state space $\Xi$, it is easy to check that the
infinitesimal generator of the Markov process with the absorption state
$\Delta$ is given by%
\begin{equation}
\Theta=\left(
\begin{array}
[c]{cc}%
0 & \mathbf{0}\\
T^{0} & T
\end{array}
\right)  , \label{equ111}%
\end{equation}
where $T\mathbf{e}+T^{0}=0$,%
\[
T^{0}=\left(
\begin{array}
[c]{c}%
T_{k}^{0}\\
T_{k-1}^{0}\\
\vdots\\
T_{3}^{0}\\
T_{2}^{0}\\
T_{1}^{0}%
\end{array}
\right)  ,T=\left(
\begin{array}
[c]{cccccc}%
T_{k,k} & T_{k,k-1} &  &  &  & \\
& T_{k-1,k-1} & T_{k-1,k-2} &  &  & \\
&  & \ddots & \ddots &  & \\
&  &  & T_{3,3} & T_{3,2} & \\
&  &  &  & T_{2,2} & T_{2,1}\\
&  &  &  &  & T_{1,1}%
\end{array}
\right)  ,
\]%
\[
T_{i}^{0}=\theta_{1}I,i=k,\ldots,2,\text{ \ }T_{1}^{0}=D_{2}+\theta_{1}I,
\]
for $j=k,k-1\ldots,2,1,$%
\[
T_{j,j}=C_{2}-j\theta_{1}I
\]
and for $i=k,k-1\ldots,3,2,$%
\[
T_{i,i-1}=D_{2}+\left(  i-1\right)  \theta_{1}I\mathbf{.}%
\]

Let $\left(  \alpha_{\Delta},\overrightarrow{\alpha}\right)  $ denote the
initial probability vector of the Markov process $\Theta$ with the absorption
state $\Delta$, $\alpha_{\Delta}=0$, $\overrightarrow{\alpha}=\left(
\alpha_{2},\mathbf{0},\ldots,\mathbf{0},\mathbf{0}\right)  $, where
$\alpha_{2}$ is the stationary probability vector of the Markov process
$C_{2}+D_{2}$. Note that the initial probability vector $\left(
\alpha_{\Delta},\overrightarrow{\alpha}\right)  $ shows that the Markov
process $\Theta$ begins at state $\left(  k,k-1\right)  $ with the phase
probability vector $\alpha_{2}$ of the MAP with irreducible matrix
representation $\left(  C_{2},D_{2}\right)  $.

The following theorem uses the phase-type distribution to provide expression
for the probability distribution of the sojourn time $W_{A}^{\left(  k\right)
}$.

\begin{The}
\label{The:PH} The probability distribution of the sojourn time $W_{A}%
^{\left(  k\right)  }$ is of phase-type with an irreducible representation
$\left(  \overrightarrow{\alpha},T\right)  $ of order $km_{2}$, and
\[
F\left(  t\right)  =P\left\{  W_{A}^{\left(  k\right)  }\leq t\right\}
=1-\overrightarrow{\alpha}\exp\left\{  Tt\right\}  \mathbf{e,}\text{ \ }%
t\geq0.
\]
Also, the average sojourn time $W_{A}^{\left(  k\right)  }$ is given by%
\[
E\left[  W_{A}^{\left(  k\right)  }\right]  =-\overrightarrow{\alpha}%
T^{-1}\mathbf{e}.
\]
\end{The}

\textbf{Proof.}\textit{ }For $\left(  m,n\right)  \in\left\{  \left(
k-i,k-i-1\right)  :i=0,1,\ldots,k-1\right\}  $ and $j_{2}\in\left\{
1,2,\ldots,m_{2}\right\}  $, we write%
\[
q_{m,n;j_{2}}\left(  t\right)  =P\left\{  M\left(  t\right)  =m,N_{1}\left(
t\right)  =n;\text{ }J_{2}\left(  t\right)  =j_{2}\right\}  ,
\]
which is the state probability that the Markov process $\Theta$ with the
absorbing state $\Delta$ is at state $\left(  m,n;j_{2}\right)  $ at time
$t\geq0$ before absorbed to state $\Delta$. We write%
\[
\mathbf{q}_{m,n}\left(  t\right)  =\left(  q_{m,n;1}\left(  t\right)
,q_{m,n;2}\left(  t\right)  ,\ldots,q_{m,n;m_{2}}\left(  t\right)  \right)  ,
\]%
\[
\overrightarrow{\mathbf{q}}\left(  t\right)  =\left(  \mathbf{q}%
_{k,k-1}\left(  t\right)  ,\mathbf{q}_{k-1,k-2}\left(  t\right)
,\ldots,\mathbf{q}_{2,1}\left(  t\right)  ,\mathbf{q}_{1,0}\left(  t\right)
\right)  .
\]
By using the Chapman-Kolmogorov forward differential equation, we obtain
\begin{equation}
\frac{\text{d}}{\text{d}t}\overrightarrow{\mathbf{q}}\left(  t\right)
=\overrightarrow{\mathbf{q}}\left(  t\right)  T \label{equ20}%
\end{equation}
with the initial condition%
\begin{equation}
\overrightarrow{\mathbf{q}}\left(  0\right)  =\overrightarrow{\alpha}.
\label{equ21}%
\end{equation}
By using $\alpha_{\Delta}=0$, it follows from (\ref{equ20}) and (\ref{equ21})
that
\begin{equation}
\overrightarrow{\mathbf{q}}\left(  t\right)  =\overrightarrow{\alpha}e^{Tt}.
\label{equ22}%
\end{equation}
Note that $\overrightarrow{\mathbf{q}}\left(  0\right)  \mathbf{e}=1$, this
gives \qquad%
\[
F\left(  t\right)  =P\left\{  W_{A}^{\left(  k\right)  }\leq t\right\}
=1-\overrightarrow{\alpha}\exp\left\{  Tt\right\}  \mathbf{e,}\text{ \ }%
t\geq0.
\]
It is easy to check that%
\[
E\left[  W_{A}^{\left(  k\right)  }\right]  =-\overrightarrow{\alpha}%
T^{-1}\mathbf{e}.
\]
This completes the proof. $\blacksquare$

The following lemma is useful for computing the inverse matrix $T^{-1}$.

\begin{Lem}
\label{Lem:Inv}The matrix $T_{j,j}=C_{2}-j\theta_{1}I$ is invertible for each
$j=k,k-1,\ldots,2,1$.
\end{Lem}

\textbf{Proof.}\textit{ }For the MAP with irreducible matrix representation
$\left(  C_{2},D_{2}\right)  $ of order $m_{2}$, note that $D_{2}\gvertneqq0$,
each diagonal element of $C_{2}$ is negative while its nondiagonal elements
are nonnegative, $C_{2}\mathbf{e}\lneqq0$, and $\left(  C_{2}+D_{2}\right)
\mathbf{e}=0$, thus it is easy to see that $C_{2}$ is a diagonally dominant
$M$-matrix. Let $\xi_{1},\xi_{2},\ldots,\xi_{m_{2}}$ be the $m_{2}$
eigenvalues of the matrix $C_{2}$, and $\operatorname{Re}\left(  \xi
_{r}\right)  $ the real part of the eigenvalue $\xi_{r}$ for $1\leq r\leq
m_{2}$. From the theory of diagonally dominant $M$-matrices, it is well-known
that $\operatorname{Re}\left(  \xi_{r}\right)  <0$ for $1\leq r\leq m_{2}$.
Since the $m_{2}$ eigenvalues of the matrix $T_{j,j}=C_{2}-j\theta_{1}I$ are
given by%
\[
\xi_{1}-j\theta_{1},\xi_{2}-j\theta_{1},\ldots,\xi_{m_{2}}-j\theta_{1},
\]
this gives that for $1\leq r\leq m_{2}$, $\operatorname{Re}\left(  \xi
_{r}-j\theta_{1}\right)  =\operatorname{Re}\left(  \xi_{r}\right)
-j\theta_{1}<0$ so that $\xi_{r}-j\theta_{1}\neq0$. Since%
\[
\det\left(  T_{j,j}\right)  =\prod_{r=1}^{m_{2}}\left(  \xi_{r}-j\theta
_{1}\right)  \neq0,
\]
this shows that the matrix $T_{j,j}=C_{2}-j\theta_{1}I$ is invertible for each
$j=k,k-1,\ldots,2,1$. This completes the proof. $\blacksquare$

To compute the average sojourn time $E\left[  W_{A}^{\left(  k\right)
}\right]  $, it is a key to deal with the inverse matrix $T^{-1}$. Let%
\[
T^{-1}=\left(
\begin{array}
[c]{ccccc}%
X_{k,k} & X_{k,k-1} & X_{k,k-2} & \cdots & X_{k,1}\\
& X_{k-1,k-1} & X_{k-1,k-2} & \cdots & X_{k-1,1}\\
&  & X_{k-2,k-2} & \cdots & X_{k-2,1}\\
&  &  & \ddots & \vdots\\
&  &  &  & X_{1,1}%
\end{array}
\right)  ,
\]
Then by using $T\cdot T^{-1}=I$ and Lemma \ref{Lem:Inv}, we obtain that for
$i=k,k-1,\ldots,2,1,$%
\[
X_{i,i}=T_{i,i}^{-1},
\]%
\[
X_{i,i-1}=-T_{i,i}^{-1}T_{i,i-1}T_{i-1,i-1}^{-1},
\]
and for $j=1,2,\ldots,i-1,$%
\[
X_{i,i-j}=\left(  -1\right)  ^{j}T_{i,i}^{-1}T_{i,i-1}T_{i-1,i-1}%
^{-1}T_{i-1,i-2}\cdots T_{i-j+1,i-j}T_{i-j,i-j}^{-1},
\]
where%
\[
T_{j,j}^{-1}=\left(  C_{2}-j\theta_{1}I\right)  ^{-1}.
\]

The following corollary provides expression for the average sojourn time
$W_{A}^{\left(  k\right)  }$.

\begin{Cor}%
\[
E\left[  W_{A}^{\left(  k\right)  }\right]  =\alpha_{2}\sum\limits_{j=0}%
^{k-1}\left(  -1\right)  ^{j+1}T_{k,k}^{-1}T_{k,k-1}T_{k-1,k-1}^{-1}\cdots
T_{k-j+1,k-j}T_{k-j,k-j}^{-1}\mathbf{e.}%
\]
\end{Cor}

\textbf{Proof.}\textit{ }It is easy to see Theorem \ref{The:PH} that%
\begin{align*}
E\left[  W_{A}^{\left(  k\right)  }\right]   &  =-\overrightarrow{\alpha
}T^{-1}\mathbf{e}\\
&  =-\left(  \alpha_{2},\mathbf{0},\ldots,\mathbf{0},\mathbf{0}\right)
\left(
\begin{array}
[c]{ccccc}%
X_{k,k} & X_{k,k-1} & X_{k,k-2} & \cdots & X_{k,1}\\
& X_{k-1,k-1} & X_{k-1,k-2} & \cdots & X_{k-1,1}\\
&  & X_{k-2,k-2} & \cdots & X_{k-2,1}\\
&  &  & \ddots & \vdots\\
&  &  &  & X_{1,1}%
\end{array}
\right)  \left(
\begin{array}
[c]{c}%
\mathbf{e}\\
\mathbf{e}\\
\mathbf{e}\\
\vdots\\
\mathbf{e}%
\end{array}
\right) \\
&  =-\alpha_{2}\sum\limits_{j=0}^{k-1}X_{k,k-j}\mathbf{e}\\
&  =\alpha_{2}\sum\limits_{j=0}^{k-1}\left(  -1\right)  ^{j+1}T_{k,k}%
^{-1}T_{k,k-1}T_{k-1,k-1}^{-1}\cdots T_{k-j+1,k-j}T_{k-j,k-j}^{-1}\mathbf{e.}%
\end{align*}
This completes the proof. $\blacksquare$

Note that $N_{1}\left(  t\right)  $ and $N_{2}\left(  t\right)  $ are the
numbers of A- and B-customers in the double-ended queue at time $t\geq0$,
respectively, and $N\left(  t\right)  =N_{1}\left(  t\right)  -N_{2}\left(
t\right)  $. Since the double-ended queue with impatient customers is always
irreducible and positive recurrent, we write%
\[
\mathbf{N}_{1}=\lim_{t\rightarrow+\infty}N_{1}\left(  t\right)  ,\text{
}\mathbf{N}_{2}=\lim_{t\rightarrow+\infty}N_{2}\left(  t\right)  ,\text{
}\mathbf{N}=\lim_{t\rightarrow+\infty}N\left(  t\right)  ,\text{ a.s..}%
\]
Let%
\[
\gamma_{k}=P\left\{  \mathbf{N}=k\right\}  ,\text{ }k=\ldots
,-2,-1,0,1,2,\ldots.
\]
Since $\mathbf{\pi}=\left(  \ldots,\mathbf{\pi}_{-2},\mathbf{\pi}%
_{-1},\mathbf{\pi}_{0},\mathbf{\pi}_{1},\mathbf{\pi}_{2},\ldots\right)  $ is
the stationary probability vector of the bilateral QBD process $Q$, we have%
\[
\gamma_{k}=\mathbf{\pi}_{k}\mathbf{e},\text{ }k=\ldots,-2,-1,0,1,2,\ldots.
\]

The following theorem provides expression for the average sojourn time of any
arriving A-customer in the double-ended queue.

\begin{The}
Note that the double-ended queue with impatient customers is always stable,
thus the average sojourn time of any arriving A-customer is given by%
\[
E\left[  W_{A}\right]  =\sum\limits_{k=1}^{\infty}\gamma_{k}E\left[
W_{A}^{\left(  k\right)  }\right]  .
\]
\end{The}

\textbf{Proof.}\textit{ }It is easy to see that%
\[
E\left[  W_{A}\right]  =\sum\limits_{k=1}^{\infty}E\left[  W_{A}|\Phi
_{k}\right]  P\left\{  \Phi_{k}\right\}  ,
\]
where for $k\geq1$%
\begin{align*}
\Phi_{k} =  &  \{\text{Any arriving A-customer finds }k-1\text{ A-customers
}\\
&  \text{(and no B-customer) in front of her at time }0\}.
\end{align*}
Specifically, we have%
\begin{align*}
\Phi_{1} =  &  \{\text{Any arriving A-customer finds neither A-customers }\\
&  \text{nor B-customer in front of her at time }0\}.
\end{align*}
Note that the double-ended queue with impatient customers is always stable, we
obtain%
\[
E\left[  W_{A}|\Phi_{k}\right]  =E\left[  W_{A}^{\left(  k\right)  }\right]
\]
and%
\[
P\left\{  \Phi_{k}\right\}  =\gamma_{k}.
\]
since%
\[
\Phi_{k}=\left\{  \text{There are }k\text{ A-customers in the double-ended
queue at time }0\right\}  .
\]
In addition, if an arriving A-customer finds no A-customer (but there is at
least one B-customer) in front of her at time $0$, it is clear that the
sojourn time of the arriving A-customer is zero because the arriving
A-customer and one B-customer can immediately match and leave the system.
Therefore, we obtain%
\[
E\left[  W_{A}\right]  =\sum\limits_{k=1}^{\infty}\gamma_{k}E\left[
W_{A}^{\left(  k\right)  }\right]  .
\]
This completes the proof. $\blacksquare$

In the remainder of this section, to understand our method for how to compute
the average sojourn times $E\left[  W_{A}^{\left(  k\right)  }\right]  $ and
$E\left[  W_{A}\right]  $, we consider a special example: The M/M/1 queue. Let
$\lambda$ and $\mu$ be the arrival and service rates, respectively, and we
assume that $\rho=\lambda/\mu<1$. Then, to compute the average sojourn times,
we have%
\[
\overrightarrow{\alpha}=\left(  1,0,0,\ldots,0\right)  ,
\]
and%
\[
T^{0}=\left(
\begin{array}
[c]{c}%
0\\
0\\
\vdots\\
0\\
\mu
\end{array}
\right)  ,T=\left(
\begin{array}
[c]{ccccc}%
-\mu & \mu &  &  & \\
& -\mu & \mu &  & \\
&  & \ddots & \ddots & \\
&  &  & -\mu & \mu\\
&  &  &  & -\mu
\end{array}
\right)  .
\]
It is easy to check that%
\[
T^{-1}=-\left(
\begin{array}
[c]{ccccc}%
\frac{1}{\mu} & \frac{1}{\mu} & \frac{1}{\mu} & \cdots & \frac{1}{\mu}\\
& \frac{1}{\mu} & \frac{1}{\mu} & \cdots & \frac{1}{\mu}\\
&  & \ddots & \ddots & \vdots\\
&  &  & \frac{1}{\mu} & \frac{1}{\mu}\\
&  &  &  & \frac{1}{\mu}%
\end{array}
\right)  .
\]
Thus we obtain%
\[
E\left[  W_{A}^{\left(  k\right)  }\right]  =-\overrightarrow{\alpha}%
T^{-1}\mathbf{e=}\frac{k}{\mu}.
\]
If $\rho<1$, then the M/M/1 queue is stable. It is well-known that%
\[
\gamma_{1}=1-\rho,\text{ \ \ }\gamma_{k}=\rho^{k-1}\left(  1-\rho\right)
,k\geq2.
\]
we have%
\begin{align*}
E\left[  W_{A}\right]   &  =\sum\limits_{k=1}^{\infty}\gamma_{k}E\left[
W_{A}^{\left(  k\right)  }\right]  \\
&  =\sum_{k=1}^{\infty}\rho^{k-1}\left(  1-\rho\right)  \frac{k}{\mu}\\
&  =\frac{1-\rho}{\mu}\sum_{k=1}^{\infty}k\rho^{k-1}\\
&  =\frac{1}{\mu\left(  1-\rho\right)  }\\
&  =\frac{1}{\mu-\lambda}.
\end{align*}

\section{Three Algorithms}

In this section, by applying the key techniques developed in Bright and Taylor
\cite{Bri:1995, Bri:1997} together with the RG-factorizations, we give three
effective algorithms for computing some stationary performance measures and
the average sojourn time $E\left[  W_{A}\right]  $.

It is worthwhile to note that the bilateral level-dependent QBD process $Q$
can be decomposed into two unilateral level-dependent QBD processes $Q_{A}$
and $Q_{B}$, as explained in Section 3. Therefore, we can use the
approximately truncated method, proposed by Bright and Taylor \cite{Bri:1995,
Bri:1997}, to determine a key truncation level $K$. In general, the truncation
level $K$ needs to be large enough such that the stationary probability of
being at all the states in or above level $K$ is sufficiently small.

For the QBD process $Q_{A}$, it follows from Bright and Taylor (1995, 1997)
that
\begin{equation}
R_{k}=\sum\limits_{l=0}^{\infty}U_{k}^{l}\prod\limits_{i=0}^{l-1}D_{k+2^{l-i}%
}^{l-1-i},\text{ }k\geq1, \label{n1}%
\end{equation}
where%
\begin{align*}
U_{k}^{0}  &  =A_{0}^{\left(  k\right)  }\left(  -A_{1}^{\left(  k+1\right)
}\right)  ^{-1},\text{ \ }k\geq1,\\
D_{k}^{0}  &  =A_{2}^{\left(  k\right)  }\left(  -A_{1}^{\left(  k-1\right)
}\right)  ^{-1},\text{ \ }k\geq1,\\
U_{k}^{l+1}  &  =U_{k}^{l}U_{k+2^{l}}^{l}\left[  I-U_{k+2^{l+1}}%
^{l}D_{k+3\cdot2^{l}}^{l}-D_{k+2^{l+1}}^{l}U_{k+2^{l}}^{l}\right]
^{-1},\text{ }l\geq0,\\
D_{k}^{l+1}  &  =D_{k}^{l}D_{k-2^{l}}^{l}\left[  I-U_{k-2^{l+1}}^{l}%
D_{k-2^{l}}^{l}-D_{k-2^{l+1}}^{l}U_{k-3\cdot2^{l}}^{l}\right]  ^{-1},\text{
}l\geq0.
\end{align*}
In addition, the $R$-measure $\left\{  R_{k}:k\geq1\right\}  $ can be
recursively computed by
\begin{equation}
R_{k}=A_{0}^{\left(  k\right)  }\left(  -A_{1}^{\left(  k+1\right)  }%
-R_{k+1}A_{2}^{\left(  k+2\right)  }\right)  ^{-1}. \label{n2}%
\end{equation}

Similarly, for the QBD process $Q_{B}$, we have%
\begin{equation}
\mathbb{R}_{k}=\sum\limits_{l=0}^{\infty}\mathbb{U}_{k}^{l}\prod
\limits_{i=0}^{l-1}\mathbb{D}_{k-2^{l-i}}^{l-1-i},\text{ }k\leq-1, \label{n3}%
\end{equation}
where%
\begin{align*}
\mathbb{U}_{k}^{0}  &  =B_{0}^{\left(  k\right)  }\left(  -B_{1}^{\left(
k-1\right)  }\right)  ^{-1},\text{ \ }k\leq-1,\\
\mathbb{D}_{k}^{0}  &  =B_{2}^{\left(  k\right)  }\left(  -B_{1}^{\left(
k+1\right)  }\right)  ^{-1},\text{ \ }k\leq-1,\\
\mathbb{U}_{k}^{l+1}  &  =\mathbb{U}_{k}^{l}\mathbb{U}_{k-2^{l}}^{l}\left[
I-\mathbb{U}_{k-2^{l+1}}^{l}\mathbb{D}_{k-3\cdot2^{l}}^{l}-\mathbb{D}%
_{k-2^{l+1}}^{l}\mathbb{U}_{k-2^{l}}^{l}\right]  ^{-1},\text{ }l\geq0,\\
\mathbb{D}_{k}^{l+1}  &  =\mathbb{D}_{k}^{l}\mathbb{D}_{k+2^{l}}^{l}\left[
I-\mathbb{U}_{k+2^{l+1}}^{l}\mathbb{D}_{k+2^{l}}^{l}-\mathbb{D}_{k+2^{l+1}%
}^{l}\mathbb{U}_{k+3\cdot2^{l}}^{l}\right]  ^{-1},\text{ }l\geq0.
\end{align*}
Also, the $R$-measure $\left\{  \mathbb{R}_{k}:k\leq-1\right\}  $ can be
recursively computed by
\begin{equation}
\mathbb{R}_{k}=B_{0}^{\left(  k\right)  }\left(  -B_{1}^{\left(  k-1\right)
}-\mathbb{R}_{k-1}B_{2}^{\left(  k-2\right)  }\right)  ^{-1}. \label{n4}%
\end{equation}

Now, our basic task is to determine such a suitable truncation level $K$. To
this end, we take a controllable accuracy $\varepsilon=10^{-20}$, and choose a
sequence $\left\{  \zeta_{k}:k\geq1\right\}  $ of positive integers with
$2\leq\zeta_{0}<\zeta_{1}<\zeta_{2}<\cdots$, for example, $\zeta_{k}=10\left(
k+1\right)  $ for $k\geq0$.

To determine the truncation level $K$, we first design Algorithm one to give
an iterative approximation by means of (\ref{n1}), (\ref{n2}), (\ref{n3}) and
(\ref{n4}). Then we design Algorithm two to compute those stationary
performance measures given in the end of Section 4, and Algorithm three to
compute the average sojourn time given in Section 5.

\vskip                                                      0.6cm

\textbf{Algorithm one: } Determination of a suitable truncation level $K$

\textit{Step 0: Initialization}

Set the initial value of $\mathit{n}$ as $\mathit{n}=0$.

\textit{Step 1: Determination of rate matrices $R_{K}$ and $\mathbb{R}_{-K}$}

Let $K=$ $\zeta_{n}$. Compute the rate matrix $R_{K}$ by means of (\ref{n1}),
and the rate matrix $\mathbb{R}_{-K}$ by means of (\ref{n3}).

\textit{Step 2: Determination of other rate matrices via an iterative computation}

(a) Based on the rate matrix $R_{K}$ where $K$ is determined in Step 1,
iteratively compute the rate matrices $R_{K-1},R_{K-2}$, $\ldots,R_{2},R_{1}$
by means of (\ref{n2}).

(b) Based on the rate matrix $\mathbb{R}_{-K}$ where $K$ is determined in Step
1, iteratively compute the rate matrices $\mathbb{R}_{-\left(  K-1\right)
},\mathbb{R}_{-\left(  K-2\right)  }$, $\ldots,\mathbb{R}_{-2},\mathbb{R}%
_{-1}$ by means of (\ref{n4}).

\textit{Step 3: Determination of the vectors }$\widetilde{\pi}_{-1}$\textit{,
}$\widetilde{\pi}_{0}$\textit{ and }$\widetilde{\pi}_{1}$

Based on the rate matrix $R_{1}$ and $\mathbb{R}_{-1}$ determined in Step 2,
determine the vectors $\widetilde{\pi}_{-1}$, $\widetilde{\pi}_{0}$ and
$\widetilde{\pi}_{1}$ through solving a system of linear equations
\[
\left\{
\begin{array}
[c]{l}%
\widetilde{\pi}_{-1}\mathbb{R}_{-1}B_{2}^{\left(  -2\right)  }+\widetilde{\pi
}_{-1}B_{1}^{\left(  -1\right)  }+\widetilde{\pi}_{0}B_{0}^{\left(  0\right)
}=0,\\
\widetilde{\pi}_{-1}B_{2}^{\left(  -1\right)  }+\widetilde{\pi}_{0}\left(
B_{1}^{\left(  0\right)  }+A_{1}^{\left(  0\right)  }\right)  +\widetilde{\pi
}_{1}A_{2}^{\left(  1\right)  }=0,\\
\widetilde{\pi}_{0}A_{0}^{\left(  0\right)  }+\widetilde{\pi}_{1}%
A_{1}^{\left(  1\right)  }+\widetilde{\pi}_{1}R_{1}A_{2}^{\left(  2\right)
}=0,\\
\widetilde{\pi}_{0}\mathbf{e}+\widetilde{\pi}_{1}\mathbf{e}+\widetilde{\pi
}_{-1}\mathbf{e}=1.
\end{array}
\right.
\]

\textit{Step 4: Determination of a normal constant}

Using the vectors $\widetilde{\pi}_{-1}$, $\widetilde{\pi}_{0}$ and
$\widetilde{\pi}_{1}$, the $R$-measure $\left\{  R_{l}:1\leq l\leq K\right\}
$ and $\left\{  \mathbb{R}_{l}:-K\leq l\leq-1\right\}  $, we calculate
\[
c=\frac{1}{\sum\limits_{k=-\left(  K+1\right)  }^{-2}\widetilde{\pi}%
_{-1}\mathbb{R}_{-1}\mathbb{R}_{-2}\cdots\mathbb{R}_{k+1}\mathbf{e}%
+\widetilde{\pi}_{-1}\mathbf{e}+\widetilde{\pi}_{0}\mathbf{e}+\widetilde{\pi
}_{1}\mathbf{e}+\sum_{k=2}^{K+1}\widetilde{\pi}_{1}R_{1}R_{2}\cdots
R_{k-1}\mathbf{e}}.
\]

\textit{Step 5: Checking convergence}

If there exists a positive integer $K=$ $\zeta_{n}$ such that%
\[
\pi_{-K-1}+\pi_{K+1}<\varepsilon,
\]
(called a stop condition), where $\pi_{-K-1}=c\widetilde{\pi}_{-1}%
\mathbb{R}_{-1}\mathbb{R}_{-2}\cdots\mathbb{R}_{-K}$ and $\pi_{K+1}%
=c\widetilde{\pi}_{1}R_{1}R_{2}\cdots R_{K}$, then $K=\zeta_{n}$, and go to
Step 6. Otherwise, let $\mathit{n}=\mathit{n}+1$, and go to Step 1.

\textit{Step 6: Output}

The algorithm stops, and the suitable truncation level is obtained as its
output, i.e., $K=\zeta_{n}$.

\vskip                 0.6cm

Once the suitable truncation level $K$ is obtained by means of Algorithm one,
we can compute some stationary performance measures given in the end of
Section 4. Here, our numerical implementations are to analyze the stationary
probabilities $P\left\{  \mathcal{Q}^{\left(  1\right)  }=0\right\}  $,
$P\left\{  \mathcal{Q}^{\left(  2\right)  }=0\right\}  $ and $P\left\{
\mathcal{Q}=0\right\}  $, and to discuss the average stationary queue lengths
$E\left[  \mathcal{Q}^{\left(  1\right)  }\right]  $, $E\left[  \mathcal{Q}%
^{\left(  2\right)  }\right]  $ and $E\left[  \mathcal{Q}\right]  $.

\textbf{Algorithm two: }Computing the stationary performance measures

\textit{Step 0: Initialization}

Using \textbf{Algorithm one}, determine the suitable truncation level $K$, the
$R$-measure $\left\{  R_{l}:1\leq l\leq K\right\}  $ and $\left\{
\mathbb{R}_{l}:-K\leq l\leq-1\right\}  $, the vectors $\widetilde{\pi}_{-1}$,
$\widetilde{\pi}_{0}$ and $\widetilde{\pi}_{1}$, and the normal constant $c$.

\textit{Step 1: Computing the stationary probability vector}

For $2\leq k\leq K+1$, we compute%
\[
\pi_{-k}=c\widetilde{\pi}_{-1}\mathbb{R}_{-1}\mathbb{R}_{-2}\cdots
\mathbb{R}_{-k+1}%
\]
and%
\[
\pi_{k}=c\widetilde{\pi}_{1}R_{1}R_{2}\cdots R_{k-1}.
\]
Also, \textit{Step 0 gives}%
\[
\pi_{-1}=c\widetilde{\pi}_{-1},\text{ }\pi_{0}=c\widetilde{\pi}_{0}%
,\text{\ }\pi_{1}=c\widetilde{\pi}_{1}.
\]

\textit{Step 2: Computing the stationary probabilities}

Using the given stationary probability vector $\left\{  \pi_{l}:-K-1\leq l\leq
K+1\right\}  $, we obtain%
\[
P\left\{  \mathcal{Q}^{\left(  1\right)  }=0\right\}  =\sum_{l=-K-1}^{0}%
\pi_{l}\mathbf{e,}%
\]%
\[
P\left\{  \mathcal{Q}^{\left(  2\right)  }=0\right\}  =\sum_{l=0}^{K+1}\pi
_{l}\mathbf{e}%
\]
and
\[
P\left\{  \mathcal{Q}=0\right\}  =\pi_{0}\mathbf{e}.
\]

\textit{Step 3: Computing the average stationary queue lengths}

Using the given stationary probability vector $\left\{  \pi_{l}:-K-1\leq l\leq
K+1\right\}  $, we obtain%
\[
E\left[  \mathcal{Q}^{\left(  1\right)  }\right]  =\sum_{l=1}^{K+1}l\pi
_{l}\mathbf{e}\text{\textbf{, }}%
\]%
\[
E\left[  \mathcal{Q}^{\left(  2\right)  }\right]  =\sum_{l=-K-1}^{-1}\left(
-l\right)  \pi_{l}\mathbf{e,}%
\]%
\[
E\left[  \mathcal{Q}\right]  =\sum_{l=1}^{K+1}l\pi_{l}\mathbf{e}\cdot
\sum_{l=1}^{K+1}\pi_{l}\mathbf{e}+\sum_{l=-K-1}^{-1}\left(  -l\right)  \pi
_{l}\mathbf{e}\cdot\sum_{l=-K-1}^{-1}\pi_{l}\mathbf{e.}%
\]

\textit{Step 4: Output}

(a) The stationary probabilities: $P\left\{  \mathcal{Q}^{\left(  1\right)
}=0\right\}  $, $P\left\{  \mathcal{Q}^{\left(  2\right)  }=0\right\}  $ and
$P\left\{  \mathcal{Q}=0\right\}  $;

(b) the average stationary queue lengths: $E\left[  \mathcal{Q}^{\left(
1\right)  }\right]  $, $E\left[  \mathcal{Q}^{\left(  2\right)  }\right]  $
and $E\left[  \mathcal{Q}\right]  $.

\vskip                 0.6cm

Once the suitable truncation level $K$ is obtained by means of Algorithm one,
we can compute the average sojourn time $E\left[  W_{A}\right]  $, given in
the Section 5.

\textbf{Algorithm three: }Computing the average sojourn time $E\left[
W_{A}\right]  $

\textit{Step 0: Initialization}

Give an initial \textit{truncation level }$K$, which is larger enough and is
given in \textbf{Algorithm one}.

\textit{Step 1: Compute the matrices: }For $j=k,k-1\ldots,2,1,$%
\[
T_{j,j}=C_{2}-j\theta_{1}I,
\]%
\[
T_{j,j}^{-1}=\left(  C_{2}-j\theta_{1}I\right)  ^{-1},
\]
and for $i=k,k-1\ldots,3,2,$%
\[
T_{i,i-1}=D_{2}+\left(  i-1\right)  \theta_{1}I\mathbf{.}%
\]
Note that $\alpha_{2}$ is the stationary probability vector of the Markov
process $C_{2}+D_{2}$, compute the average sojourn time%
\[
E\left[  W_{A}^{\left(  k\right)  }\right]  =\alpha_{2}\sum\limits_{j=0}%
^{k-1}\left(  -1\right)  ^{j+1}T_{k,k}^{-1}T_{k,k-1}T_{k-1,k-1}^{-1}\cdots
T_{k-j+1,k-j}T_{k-j,k-j}^{-1}\mathbf{e}.
\]

\textit{Step 2: Compute the stationary probabilities}%
\[
\gamma_{k}=\mathbf{\pi}_{k}\mathbf{e},\text{ \ }k=1,2,\ldots,K.
\]

\textit{Step 3: Compute the average sojourn time} $E\left[  W_{A}\right]  $%
\[
E\left[  W_{A}\right]  =\sum\limits_{k=1}^{K}\gamma_{k}E\left[  W_{A}^{\left(
k\right)  }\right]  \mathbf{.}%
\]

\textit{Step 4: Output}

The algorithm stops, and the $E\left[  W_{A}\right]  $ is obtained as its output.

\section{Numerical Examples}

In this section, by using the above three algorithms developed in Section 6,
we provide some numerical examples to show how some performance measures of
the double-ended queue are influenced by key system parameters. Also, we apply
the coupling method to give some interesting interpretations on the numerical results.

In the numerical implementations, we shall discuss three groups of interesting
issues: (1) The stationary performance measures, (2) the sojourn time
$E\left[  W_{A}\right]  $, and (3) further numerical analysis.

\vskip              0.3cm

\textbf{(a) The stationary performance measures}

By using \textbf{Algorithms one and two}, we use some numerical examples to
analyze how the stationary performance measures of the double-ended queue are
influenced by the two key system parameters $\theta_{1}$ and $\theta_{2}$,
i.e., the impatient rates of A- and B-customers, respectively.

For the two types of customers, we respectively take their MAPs with
irreducible matrix representations as follows:
\[
C_{1}=\left(
\begin{array}
[c]{cc}%
-10 & 0\\
1 & -1
\end{array}
\right)  \text{, }D_{1}=\left(
\begin{array}
[c]{cc}%
9 & 1\\
0 & 0
\end{array}
\right)  \text{;}%
\]%
\[
C_{2}=\left(
\begin{array}
[c]{cc}%
-5 & 1\\
2 & -7
\end{array}
\right)  \text{, }D_{2}=\left(
\begin{array}
[c]{cc}%
0 & 4\\
2 & 3
\end{array}
\right)  \text{.}%
\]
It is easy to compute that $\alpha_{1}=\left(  1/2,1/2\right)  $ and
$\alpha_{2}=\left(  4/9,5/9\right)  $. Hence, their stationary arrival rates
$\lambda_{1}=\alpha_{1}D_{1}\mathbf{e}=5$ and $\lambda_{2}=\alpha_{2}%
D_{2}\mathbf{e}=4\frac{5}{9}$.

\vskip                                                     0.3cm

Now, we analyze the probabilities of stationary queue lengths for the
A-customers, the B-customers and the total system, respectively.

From Figure 4, it is observed that the stationary probability $P\left\{
\mathcal{Q}^{\left(  1\right)  }=0\right\}  $ increases as $\theta_{1}$
increases, while it decreases as $\theta_{2}$ increases.

\begin{figure}[tbh]
\centering        \includegraphics[height=7cm
,width=16cm]{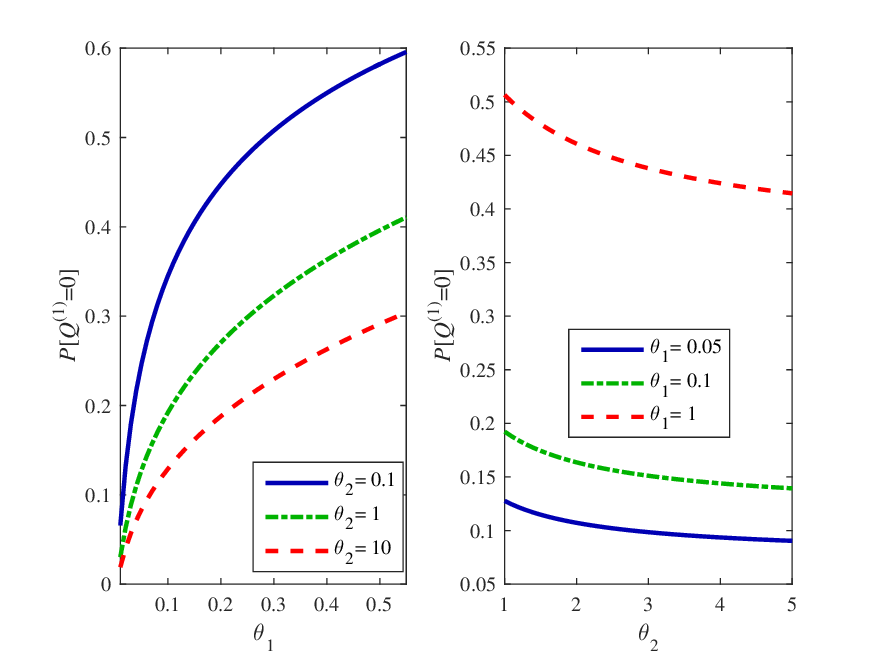}
\caption{$P\left\{  \mathcal{Q}^{\left(  1\right)  }=0\right\}  $ vs
$\theta_{1}$ and $\theta_{2}$}%
\end{figure}

From Figure 5, it is seen that the stationary probability $P\left\{
\mathcal{Q}^{\left(  2\right)  }=0\right\}  $ decreases as $\theta_{1}$
increases, while it increases as $\theta_{2}$ increases. \

\begin{figure}[tbh]
\centering        \includegraphics[height=7cm
,width=16cm]{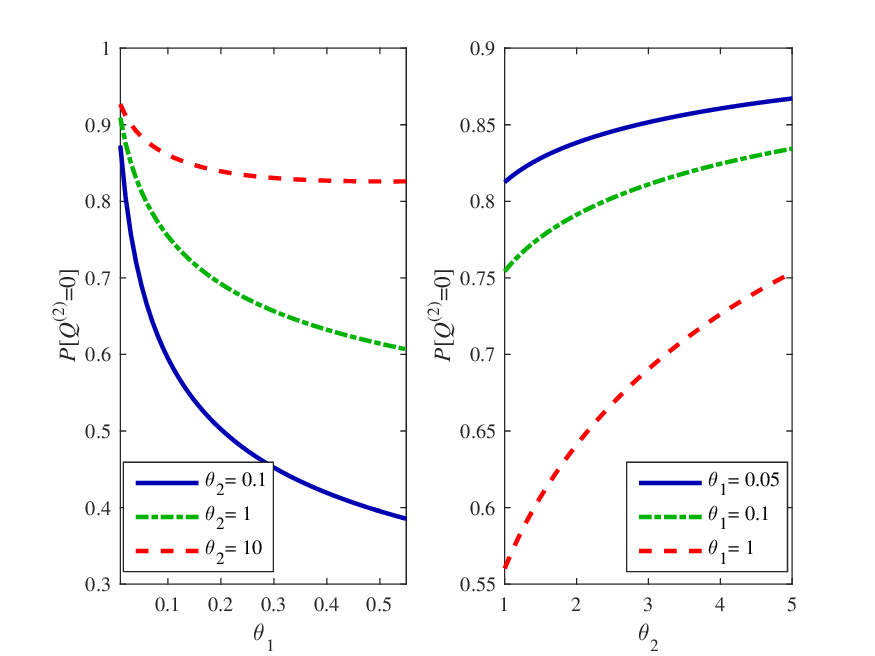}
\caption{$P\left\{  \mathcal{Q}^{\left(  2\right)  }=0\right\}  $ vs
$\theta_{1}$ and $\theta_{2}$}%
\end{figure}

\textbf{A coupling analysis: }The two numerical results can be intuitively
understood by means of the coupling method. As $\theta_{1}$ increases, more
and more A-customers are quickly leaving the system due to their impatient
behavior, thus $P\left\{  \mathcal{Q}^{\left(  1\right)  }=0\right\}  $
increases. On the other hand, as $\theta_{2}$ increases, more and more
B-customers are quickly leaving the system so that the chance that an
A-customer can match a B-customer will become smaller and smaller, hence this
leads to the decrease of $P\left\{  \mathcal{Q}^{\left(  1\right)
}=0\right\}  $.

From Figure 6, we observed an interesting phenomenon: The stationary
probability that there is neither A-customer nor B-customer, $P\left\{
\mathcal{Q}=0\right\}  $, increases as $\theta_{1}$ or $\theta_{2}$ increase.

\begin{figure}[tbh]
\centering        \includegraphics[height=7cm
,width=16cm]{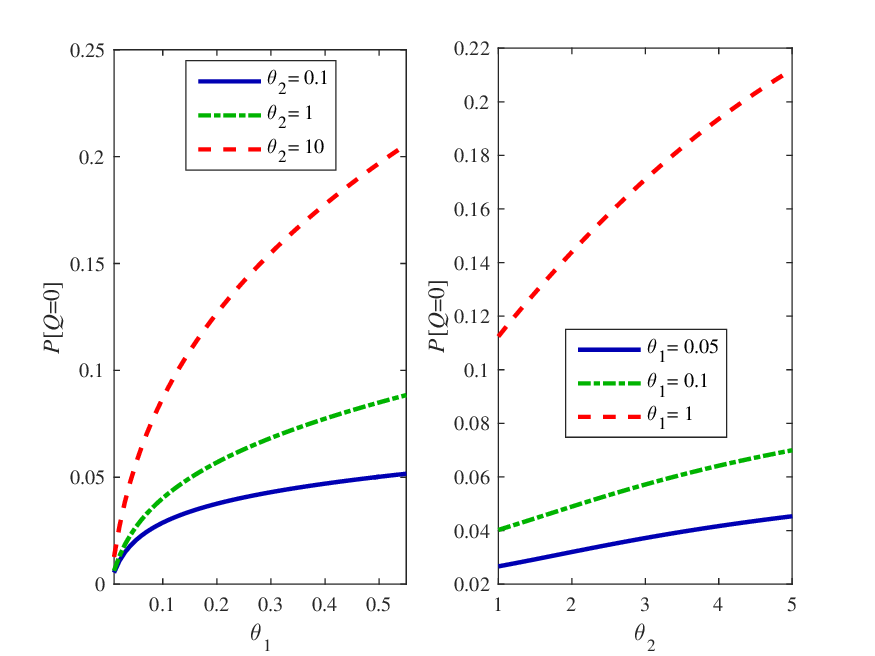}
\caption{$P\left\{  \mathcal{Q}=0\right\}  $ vs $\theta_{1}$ and $\theta_{2}$}%
\end{figure}

\vskip                 0.3cm

In what follows we discuss the average stationary queue lengths for the
A-customers, the B-customers and the total system, respectively.

Figure 7 shows that the average stationary queue length $E\left[
\mathcal{Q}^{\left(  1\right)  }\right]  $ decreases as $\theta_{1}$
increases, while it increases as $\theta_{2}$ increases. From Figure 8, we
observe that the average stationary queue length $E\left[  \mathcal{Q}%
^{\left(  2\right)  }\right]  $ increases as $\theta_{1}$ increases, while it
decreases as $\theta_{2}$ increases.

\begin{figure}[tbh]
\centering        \includegraphics[height=7cm
,width=16cm]{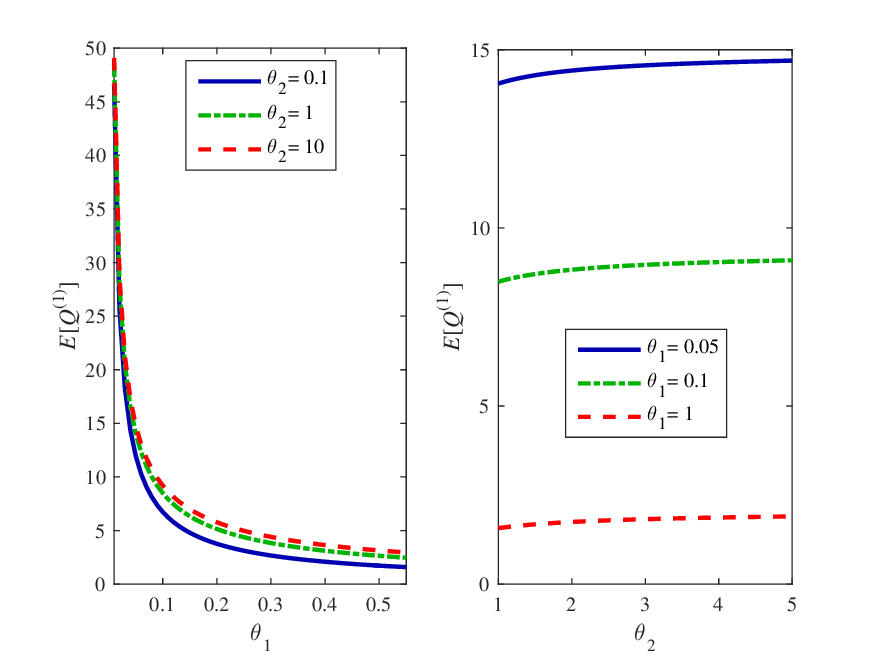}
\caption{$E\left[  \mathcal{Q}^{\left(  1\right)  }\right]  $ vs $\theta_{1}$
and $\theta_{2}$}%
\end{figure}

\begin{figure}[tbh]
\centering        \includegraphics[height=7cm
,width=16cm]{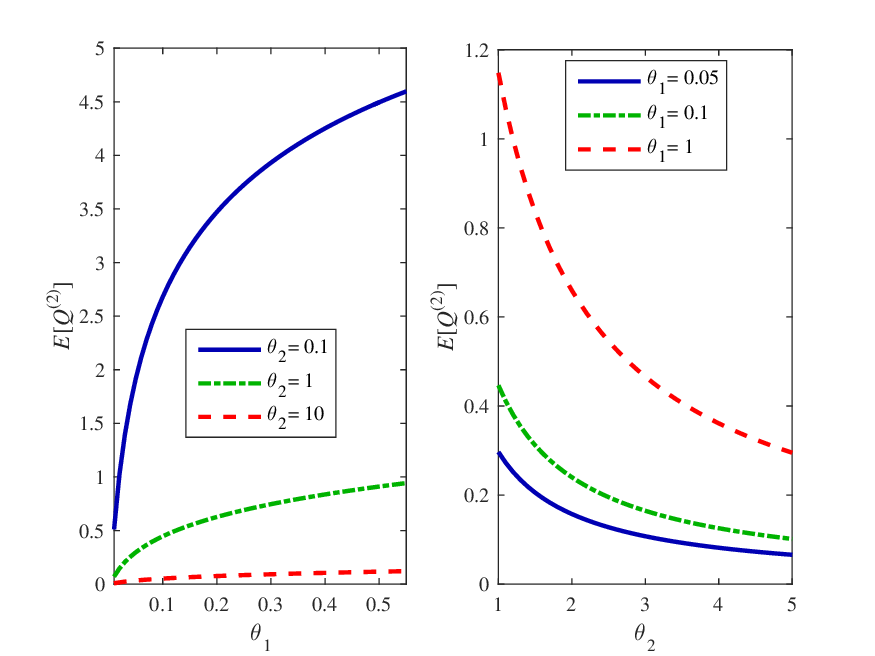}
\caption{$E\left[  \mathcal{Q}^{\left(  2\right)  }\right]  $ vs $\theta_{1}$
and $\theta_{2}$}%
\end{figure}

\textbf{A coupling analysis: }The two numerical results are intuitive. As
$\theta_{1}$ increases, more and more A-customers are quickly leaving the
system so that $E\left[  \mathcal{Q}^{\left(  1\right)  }\right]  $ decreases.
On the other hand, as $\theta_{2}$ increases, more and more B-customers are
quickly leaving the system so that the chance that an A-customer can match a
B-customer will become smaller and smaller. Thus, $E\left[  \mathcal{Q}%
^{\left(  1\right)  }\right]  $ increases as $\theta_{2}$ increases.

From Figure 9, it is seen that the average stationary queue length $E\left[
\mathcal{Q}\right]  $ decreases as $\theta_{1}$ increases. But, $E\left[
\mathcal{Q}\right]  $ has a more consistent behavior as $\theta_{2}$ increases.

\begin{figure}[tbh]
\centering        \includegraphics[height=7cm
,width=16cm]{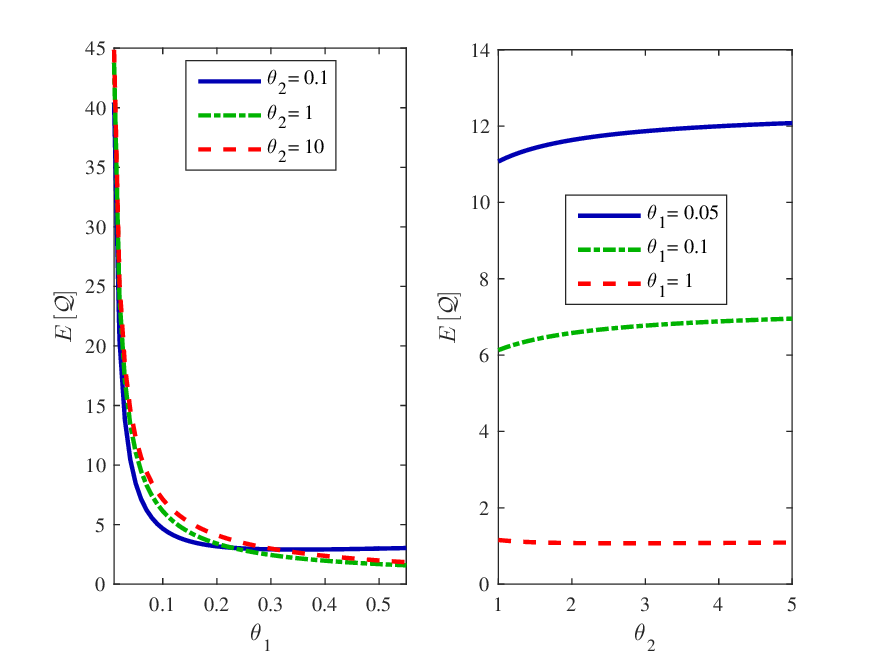}
\caption{$E\left[  \mathcal{Q}\right]  $ vs $\theta_{1}$ and $\theta_{2}$}%
\end{figure}

\vskip              0.3cm

\textbf{(b) The sojourn time}

In this subsection, by using \textbf{Algorithm three}, we use some numerical
examples to indicate how the average sojourn time $E\left[  W_{A}\right]  $
depends on the impatient rate $\theta_{1}$ and the number $k$ of A-customers
in the double-ended queue.

For the two types of customers, we respectively take their MAPs with
irreducible matrix representations as follows:
\[
C_{1}=\left(
\begin{array}
[c]{cc}%
-10 & 0\\
1 & -1
\end{array}
\right)  ,\text{ }D_{1}=\left(
\begin{array}
[c]{cc}%
9 & 1\\
0 & 0
\end{array}
\right)  ;
\]%
\[
C_{2}=\left(
\begin{array}
[c]{cc}%
-5 & 1\\
2 & -7
\end{array}
\right)  ,\text{ }D_{2}=\left(
\begin{array}
[c]{cc}%
0 & 4\\
2 & 3
\end{array}
\right)  .
\]
It is easy to compute that $\alpha_{1}=\left(  1/2,1/2\right)  $ and
$\alpha_{2}=\left(  4/9,5/9\right)  $, hence their same stationary arrival
rates are given by $\lambda_{1}=\alpha_{1}D_{1}\mathbf{e}=5$ and $\lambda
_{2}=\alpha_{2}D_{2}\mathbf{e}=4\frac{5}{9}$.

Figure 10 indicates that when $\theta_{2}=0.1,1,10$, $E\left[  W_{A}\right]  $
decreases as $\theta_{1}$ increases for $\theta_{1}\in\left[
0.01,0.55\right]  $. Figure 11 shows that when $\theta_{1}=1,2,5$, $E\left[
W_{A}^{\left(  k\right)  }\right]  $ increases as $k$ increases for
$k=5,6,\ldots,50$. As $\theta_{1}$ increases, the arriving A-customer can
quickly leave the system due to her impatient behavior, thus $E\left[
W_{A}\right]  $ decreases. As $k$ increases, the total matching time length of
the $k-1$ A-customer in front of this arriving A-customer will increase. This
shows that $E\left[  W_{A}^{\left(  k\right)  }\right]  $ increases, as $k$ increases.

\begin{figure}[tbh]
\centering
\includegraphics[height=6cm
,width=10cm]{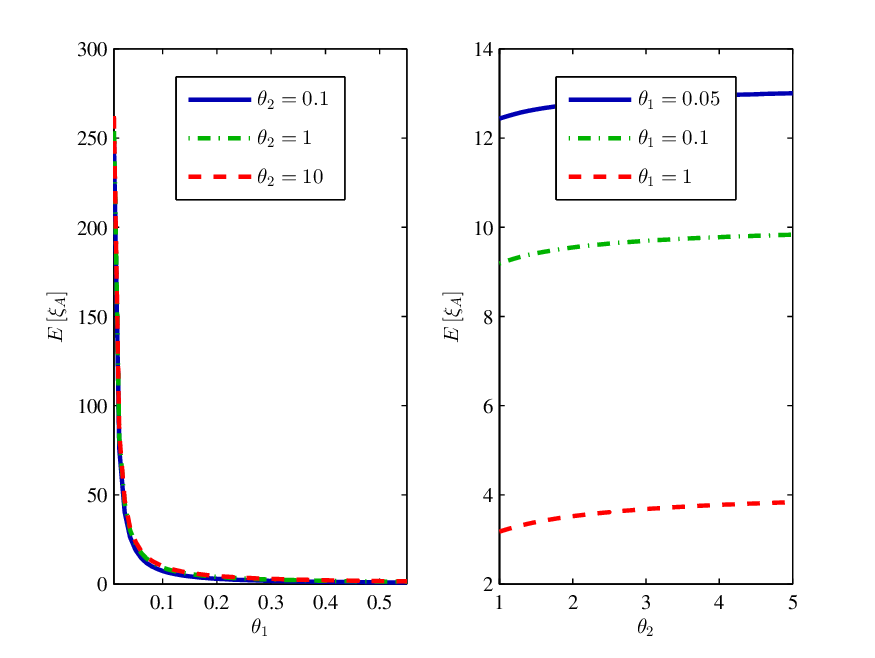} \caption{$E\left[
W_{A}\right]  $ vs $\theta_{1}$}%
\end{figure}

\begin{figure}[tbh]
\centering
\includegraphics[height=6cm
,width=10cm]{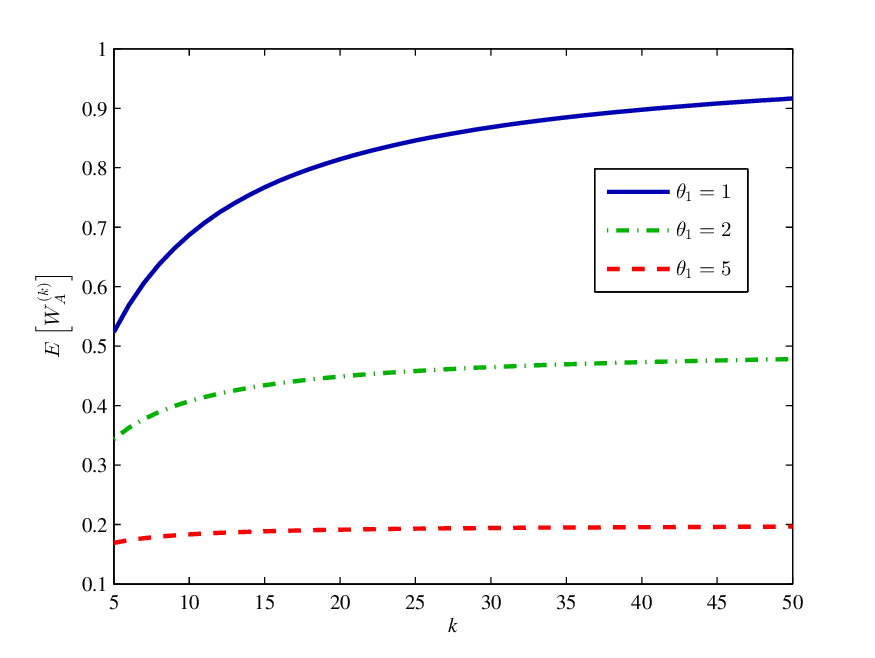} \caption{$E\left[
W_{A}^{\left(  k\right)  }\right]  $ vs $k$}%
\end{figure}

Now, we consider another example with the two MAPs of order 4 (i.e., the order
of the irreducible matrix representation is $4$), both of which have the same
stationary arrival rates as that in the previous example, i.e., $\lambda
_{1}=\alpha_{1}D_{1}\mathbf{e}=5$ and $\lambda_{2}=\alpha_{2}D_{2}%
\mathbf{e}=4\frac{5}{9}$. Also, the two MAPs of order 4 have the irreducible
matrix representations as follows:%
\[
C_{1}=\left(
\begin{array}
[c]{cccc}%
-7 & 0 & 2 & 0\\
2 & -7 & 3 & 0\\
0 & 0 & -10 & 0\\
2 & 1 & 2 & -8
\end{array}
\right)  ,\text{ }D_{1}=\left(
\begin{array}
[c]{cccc}%
0 & 5 & 0 & 0\\
0 & 1 & 1 & 0\\
0 & 0 & 2 & 8\\
3 & 0 & 0 & 0
\end{array}
\right)  ;
\]%
\[
C_{2}=\left(
\begin{array}
[c]{cccc}%
-2 & 0 & 0 & 0\\
0 & -7 & 0 & 0\\
0 & 0 & -15 & 0\\
1/2 & 0 & 5/2 & -5
\end{array}
\right)  ,\text{ }D_{2}=\left(
\begin{array}
[c]{cccc}%
0 & 2 & 0 & 0\\
0 & 3 & 4 & 0\\
3 & 0 & 2 & 10\\
2 & 0 & 0 & 0
\end{array}
\right)  .
\]
It is easy to check that $\alpha_{1}=\left(  1/4,1/4,1/4,1/4\right)  $ and
$\alpha_{2}=\left(  4/9,2/9,1/9,2/9\right)  $. Hence we obtain that
$\lambda_{1}=\alpha_{1}D_{1}\mathbf{e}=5$ and $\lambda_{2}=\alpha_{2}%
D_{2}\mathbf{e}=4\frac{5}{9}$.

\begin{figure}[tbh]
\centering
\includegraphics[height=6cm
,width=10cm]{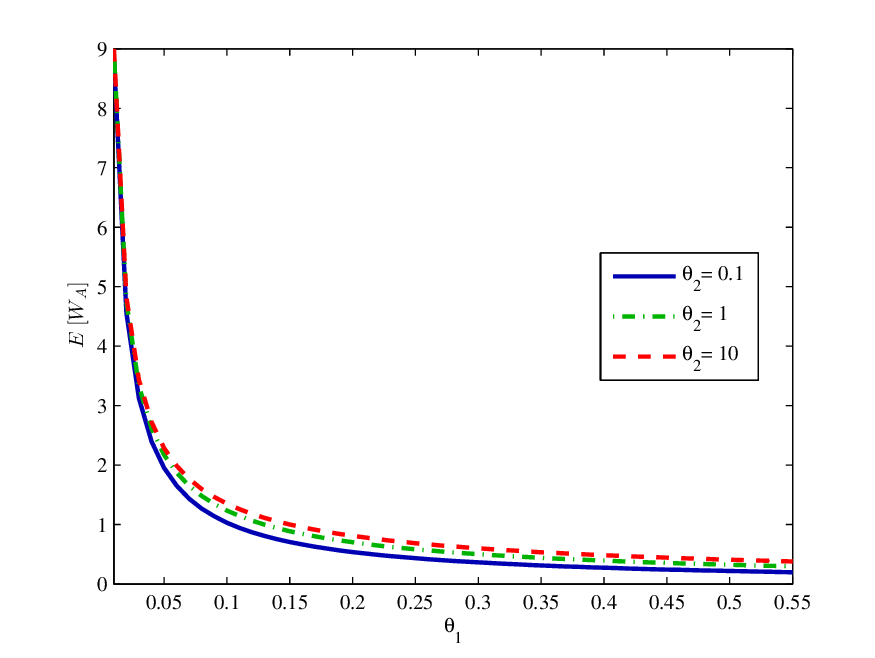} \caption{$E\left[
W_{A}\right]  $ vs $\theta_{1}$}%
\end{figure}

\begin{figure}[tbh]
\centering
\includegraphics[height=6cm
,width=10cm]{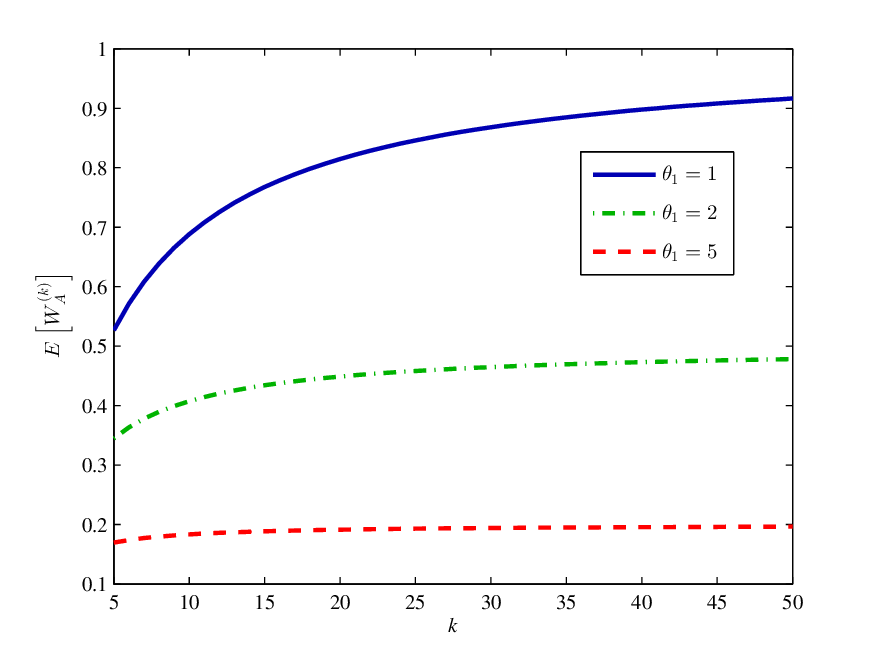} \caption{$E\left[
W_{A}^{\left(  k\right)  }\right]  $ vs $k$}%
\end{figure}

It is observed from Figure 12 that when $\theta_{2}=0.1,1,10$, $E\left[
W_{A}\right]  $ decreases as $\theta_{1}$ increases for $\theta_{1}\in\left[
0.01,0.55\right]  $. Figure 13 shows that when $\theta_{1}=1,2,5$, $E\left[
W_{A}^{\left(  k\right)  }\right]  $ increases as $k$ increases for
$k=5,6,\ldots,50$. Thus, it is easy to see that there are two identical
changing trends of $E\left[  W_{A}^{\left(  k\right)  }\right]  $ through
observing Figures 10 and 12, and Figures 11 and 13, respectively.

\vskip              0.3cm

\textbf{(c) Further numerical analysis}

Now, we show that the above three algorithms developed in Section 6 are
effective in numerical analysis of the double-ended queue when the two MAPs
have different orders. To do this, we discuss three different cases: (1) The
Poisson processes, (2) the MAP of order $2$, and (3) the MAP of order $4$.
Also, it is seen from our numerical experiments that we can easily deal with
the cases with MAPs of higher order.

We consider the double-ended queue from two groups of different impatient
rates: $\left(  \theta_{1},\theta_{2}\right)  =\left(  0.25,1\right)  $, and
$\left(  \theta_{1},\theta_{2}\right)  =\left(  0.75,1\right)  $,
respectively. While for the arrival processes in two sides, we consider three
different cases with the same stationary arrival rates: $\lambda_{1}%
=\alpha_{1}D_{1}\mathbf{e}=5$ and $\lambda_{2}=\alpha_{2}D_{2}\mathbf{e}%
=4\frac{5}{9}$. Further, we take the two arrival processes as follows:

(1) Two Poisson processes with arrival rates $\lambda_{1}=5$ and $\lambda
_{2}=4\frac{5}{9}$, respectively.

(2) Two MAPs of order 2 whose irreducible matrix representations are given by
\[
C_{1}=\left(
\begin{array}
[c]{cc}%
-10 & 0\\
1 & -1
\end{array}
\right)  ,\text{ }D_{1}=\left(
\begin{array}
[c]{cc}%
9 & 1\\
0 & 0
\end{array}
\right)  ;
\]%
\[
C_{2}=\left(
\begin{array}
[c]{cc}%
-5 & 1\\
2 & -7
\end{array}
\right)  ,\text{ }D_{2}=\left(
\begin{array}
[c]{cc}%
0 & 4\\
2 & 3
\end{array}
\right)  ,
\]
respectively.

(3) Two MAP of order 4 whose irreducible matrix representations are given by%
\[
C_{1}=\left(
\begin{array}
[c]{cccc}%
-7 & 0 & 2 & 0\\
2 & -7 & 3 & 0\\
0 & 0 & -10 & 0\\
2 & 1 & 2 & -8
\end{array}
\right)  ,\text{ }D_{1}=\left(
\begin{array}
[c]{cccc}%
0 & 5 & 0 & 0\\
0 & 1 & 1 & 0\\
0 & 0 & 2 & 8\\
3 & 0 & 0 & 0
\end{array}
\right)  ;
\]%
\[
C_{2}=\left(
\begin{array}
[c]{cccc}%
-2 & 0 & 0 & 0\\
0 & -7 & 0 & 0\\
0 & 0 & -15 & 0\\
1/2 & 0 & 5/2 & -5
\end{array}
\right)  ,\text{ }D_{2}=\left(
\begin{array}
[c]{cccc}%
0 & 2 & 0 & 0\\
0 & 3 & 4 & 0\\
3 & 0 & 2 & 10\\
2 & 0 & 0 & 0
\end{array}
\right)  ,
\]
respectively.

Based on the three different arrival processes, Table 1 provides a numerical
comparison for the performance measures of the double-ended queue.

\begin{center}
\textbf{Table 1}~~Comparison of performance measures under different arrival
processes\newline \resizebox{\textwidth}{20mm}{
\begin{tabular}%
{ccccccccc}
\hline
$\left(  \theta_{1},\,\theta_{2}%
\right)  $&Arrival\,Process & $P\left\{  \mathcal{Q}^{\left(  1\right)
}%
=0\right\}  $ & $P\left\{  \mathcal{Q}^{\left(  2\right)  }%
=0\right\}  $ &
$P\left\{  \mathcal{Q}=0\right\}  $ & $E\left[  \mathcal{Q}%
^{\left(
1\right)  }\right]  $ & $E\left[  \mathcal{Q}^{\left(  2\right)  }%
\right]  $ &
$E\left[  \mathcal{Q}\right]  $& $E\left[  W_{A}%
\right]  $\\
\hline
\multirow{3}*{(0.25,\,1)}%
&Poisson & 0.2850 & 0.8174 &
0.1024 & 3.3181 & 0.3851 & 2.4429& 0.5423\\
&MAP of order 2 & 0.3119 & 0.6600 & 0.0526& 4.1073 & 0.7565 & 2.6912 &
0.5995\\
&MAP of order 4 & 0.2853 & 0.8068 & 0.0966 & 3.4584 & 0.4222 & 2.5357&
0.5827\\
\hline
\multirow{3}%
*{(0.75,\,1)}%
&Poisson & 0.4699 & 0.6989 &
0.1688 & 1.4392 & 0.6350 & 0.9542& 0.2016\\
&MAP of size 2 & 0.5103 & 0.5338 & 0.0861& 1.6470 & 1.2327 & 1.2601 &
0.2038\\
&MAP of size 4 & 0.4646 & 0.6890 & 0.1571 & 1.5095 & 0.6892 & 1.0149&
0.2243\\ \hline
\end{tabular}%
}
\end{center}

From Table 1, it is seen that the performance measures with two Poisson
processes are close to that with two MAPs of order 4; while they are
significantly different from that with two MAPs of order 2.

From the numerical computation in Table 1, it is seen that our three
algorithms (by using the matrix-analytic method and the RG-factorizations),
developed in Section 6, can effectively deal with the MAPs of higher orders in
numerical analysis of the double-ended queue, for example, the MAP of order
10, the MAP of order 20, and so on.

\vskip             0.3cm

In what follows we further give a numerical example to compare the results of
this paper with that given by the multi-layer MMFF method given in Wu and He
\cite{WuH:2020} and by the diffusion method given in Liu et al.
\cite{Liu:2014}. We assume that the interarrival times of the arrival
processes intwo sides are taken as follows:

Case one: The two exponential distributions with parameters $\lambda_{1}=1$
and $\lambda_{2}=2$, respectively.

Case two: The two Erlang distributions of order $2$ whose parameters are given
by $2\lambda_{1}=2$ and $2\lambda_{2}=4$, respectively.

If the two interarrival times are exponential with parameter $\lambda_{1}=1$
and $\lambda_{2}=2$ respectively, then both of them correspond to the MAPs of
order $1$, that is,%
\begin{align*}
C_{1}  &  =-1,\text{ }D_{1}=1;\text{ }\\
C_{2}  &  =-2,\text{ }D_{2}=2\text{.}%
\end{align*}

If the two Erlang distributions of order $2$ whose parameters are given by
$2\lambda_{1}=2$ and $2\lambda_{2}=4$ respectively, then both of them
correspond to the MAPs of order $2$, that is,%
\[
C_{1}=\left(
\begin{array}
[c]{cc}%
-2 & 2\\
0 & -2
\end{array}
\right)  ,\text{ }D_{1}=\left(
\begin{array}
[c]{cc}%
0 & 0\\
2 & 0
\end{array}
\right)  ;
\]%
\[
C_{2}=\left(
\begin{array}
[c]{cc}%
-4 & 4\\
0 & -4
\end{array}
\right)  ,\text{ }D_{2}=\left(
\begin{array}
[c]{cc}%
0 & 0\\
4 & 0
\end{array}
\right)  .
\]
Note that $\alpha_{1}=\alpha_{2}=\left(  1/2,1/2\right)  $, thus we obtain
that $\lambda_{1}=\alpha_{1}D_{1}\mathbf{e}=1$ and $\lambda_{2}=\alpha
_{2}D_{2}\mathbf{e}=2$.

\vskip     0.2cm

\begin{center}
\textbf{Table 2}~~Numerical comparison of average stationary queue length
difference $E\left[  \mathcal{Q}^{\left(  1\right)  }\right]  -E\left[
\mathcal{Q}^{\left(  2\right)  }\right]  $ among three different methods
\newline \resizebox{\textwidth}{34mm}  {
\begin{tabular}
[c]{ccccccccc}%
\hline
\multicolumn{2}{c}{Parameters} & \multirow{2}  *{Our bilateral QBD}%

&
\multicolumn{2}{c}{Multi-layer MMFF method of Wu and He (2021)}%

&
\multicolumn{4}{c}{Diffusion method of Liu et al. (2014)}%
\\
\cmidrule  (r){1-2}\cmidrule  (r){4-5}\cmidrule  (r){6-9}%
$\left(  \lambda
_{1},\,\lambda_{2}%
\right)  \left(  1,\,2\right)  $ & $\left(  \theta
_{1},\,\theta_{2}%
\right)  $ &  & $M=N=1000$ & $M=N=2000$ & Simulation &
Poisson & Diffusion \thinspace1 & Diffusion \thinspace2\\\hline
Erlang(2) & (1,\thinspace2) & -0.4491 & -0.4396 & -0.4334 & -0.4285 &
-0.3858 & -0.4493 & -0.5\\
&  & ($4.81\%$) & ($2.58\%$) & ($1.13\%$) & ($\pm0.0018$) & ($9.96\%$) &
($4.87\%$) & ($16.69\%$)\\
Erlang(2) & (0.1,\thinspace0.2) & -5.1495 & -5.0847 & -5.0407 & -4.9832 &
-4.9719 & -4.9983 & -5\\
&  & ($3.34\%$) & ($2.04\%$) & ($1.15\%$) & ($\pm0.0015$) & ($0.23\%$) &
($0.30\%$) & ($0.34\%$)\\
Erlang(2) & (0.01,\thinspace0.02) & -50.8400 & -50.8731 & -50.4535 & -50.089 &
-50 & -50 & -50\\
&  & ($1.50\%$) & ($1.57\%$) & ($0.73\%$) & ($\pm0.1507$) & ($0.18\%$) &
($0.18\%$) & ($0.18\%$)\\
Exponential & (1,\thinspace2) & -0.3858 & -0.4007 & -0.3933 & -0.3876 &
-0.3858 & -0.3178 & -0.5\\
&  & ($0.46\%$) & ($3.38\%$) & ($1.46\%$) & ($\pm0.002$) & ($0.46\%$) &
($18\%$) & ($29\%$)\\
Exponential & (0.1,\thinspace0.2) & -4.9719 & -5.0620 & -5.0171 & -4.9779 &
-4.9719 & -4.9776 & -5\\
&  & ($0.12\%$) & ($1.69\%$) & ($0.79\%$) & ($\pm0.0157$) & ($0.12\%$) &
($0.01\%$) & ($0.45\%$)\\
Exponential & (0.01,\thinspace0.02) & -50 & -50.8615 & -50.4406 & -49.9609 &
-50 & -50 & -50\\
&  & ($0.08\%$) & ($1.80\%$) & ($0.96\%$) & ($\pm0.142$) & ($0.08\%$) &
($0.08\%$) & ($0.08\%$)\\\hline
\end{tabular}%

}
\end{center}

In contrast with a simulation result, we define a relative error ratio of
Method A as%
\[
x\%=\frac{100\left(  \text{The numerical result given by Method A }-\text{ The
simulation result}\right)  }{\text{The simulation result}}\%\text{.}%
\]
From Columns $3$ to $9$ in Table 2, it is easy to see that the relative error
ratios (in the brackets) of our bilateral QBD process are very close to that
given by both the multi-layer MMFF method of Wu and He (2021) and the
diffusion method of Liu et al. (2014). From such a comparison as well as the
numerical computation of the well-known matrix-analytic method, our bilateral
QBD process has two advantages: (a) It can easily provide a more detailed
performance analysis of the double-ended queues, especially in the cases with
MAP (non-Poisson) inputs. (b) Our bilateral QBD process can easily provide
numerical computation of the double-ended queues with MAPs of higher order
through using the matrix-analytic method and the RG-factorizations.

\section{Concluding Remarks}

In this paper, we study a block-structured double-ended queue with two MAP
inputs and customers' impatient behavior, and show that such a double-ended
queue can be expressed as a new bilateral QBD process. Based on this finding,
we provide a detailed analysis for the block-structured double-ended queue,
including the system stability, the stationary queue length and the sojourn
time. At the same time, we develop three effective algorithms for numerically
computing performance measures of the block-structured double-ended queue,
such as the probabilities of stationary queue lengths, the average stationary
queue lengths, and the average sojourn time. Finally, we use some numerical
examples to indicate how the performance measures are influenced by key system
parameters. We believe that the methodology and results given in this paper
can be applicable to deal with more general double-ended queues in practice,
and further develop some effective algorithms for the purpose of many actual uses.

Along these lines, we will continue our future research on the following directions:

-- Consider more general double-ended queues, for example, a double-ended
queue with two BMAP inputs, a double-ended queue with a BMAP input and a
renewal-process input, and a double-ended queue with matching batch size pair
$(m, n)$.

-- Develop more bilateral block-structured Markov processes, for example,
bilateral Markov processes of GI/M/1 type, bilateral Markov processes of M/G/1
type, and so on.

-- Develop effective algorithms for analyzing bilateral block-structured
Markov processes and provide numerical analysis for more general matching
issues in practice.

-- Develop stochastic optimization and dynamic control, Markov decision
processes and stochastic game theory in the study of double-ended queues. In
this case, developing effective algorithms for dealing with optimal and
control issues of the double-ended queues.

\section*{Acknowledgements}

Quan-Lin Li was supported by the National Natural Science Foundation of China
under grants No. 71671158 and 71932002 and by Beijing Social Science
Foundation Research Base Project under grant No. 19JDGLA004.

\section*{Appendix}

This appendix provides the proofs of Theorems \ref{The1} and \ref{The2}. Our
purpose is to increase the readability of the main paper.

\textbf{(a) Proof of Theorem \ref{The1}}

It is easy to see the irreducibility of the bilateral QBD process $Q$ through
observing Figure 2, and using the two MAP inputs as well as the two
exponential impatient times with $\left(  \theta_{1},\theta_{2}\right)  >0$.

Note that the bilateral QBD process $Q$ is positive recurrent if and only if
the two unilateral QBD processes $Q_{A}$ and $Q_{B}$ are positive recurrent,
thus it is key to find some necessary and sufficient conditions for the
stability of the two unilateral QBD processes $Q_{A}$ and $Q_{B}$ by means of
Neuts' method (i.e., the mean drift technique).

For the QBD process $Q_{A}$,\textit{ }let $\mathbb{A}_{k}=A_{0}^{\left(
k\right)  }+A_{1}^{\left(  k\right)  }+A_{2}^{\left(  k\right)  }$. Then for
$k\geq1$ we obtain
\begin{align*}
\mathbb{A}_{k}  &  =I\otimes D_{1}+C_{2}\oplus C_{1}-k\theta_{1}I+D_{2}\otimes
I+k\theta_{1}I\\
&  =I\otimes D_{1}+C_{2}\oplus C_{1}+D_{2}\otimes I\\
&  =I\otimes D_{1}+I\otimes C_{1}+C_{2}\otimes I+D_{2}\otimes I\\
&  =I\otimes\left(  C_{1}+D_{1}\right)  +\left(  C_{2}+D_{2}\right)  \otimes
I\\
&  =\left(  C_{2}+D_{2}\right)  \oplus\left(  C_{1}+D_{1}\right)  ,
\end{align*}
which is independent of the positive integer $k\geq1$. Obviously,
$\mathbb{A}_{k}$ is the infinitesimal generator of the continuous-time Markov
process with $m_{1}m_{2}$ states.

Note that $\alpha_{1}$ and $\alpha_{2}$ are the stationary probability vectors
of the Markov processes $C_{1}+D_{1}$ and $C_{2}+D_{2}$, respectively. Thus,
$\alpha_{1}\left(  C_{1}+D_{1}\right)  =\mathbf{0}$, $\alpha_{1}\mathbf{e}=1$;
and $\alpha_{2}\left(  C_{2}+D_{2}\right)  =\mathbf{0}$, and $\alpha
_{2}\mathbf{e}=1$. For each $k\geq1$, we get%
\begin{align*}
\left(  \alpha_{2}\otimes\alpha_{1}\right)  \mathbb{A}_{k}  &  =\left(
\alpha_{2}\otimes\alpha_{1}\right)  \left[  \left(  C_{2}+D_{2}\right)
\oplus\left(  C_{1}+D_{1}\right)  \right] \\
&  =\left(  \alpha_{2}\otimes\alpha_{1}\right)  \left[  \left(  C_{2}%
+D_{2}\right)  \otimes I\right]  +\left(  \alpha_{2}\otimes\alpha_{1}\right)
\left[  I\otimes\left(  C_{1}+D_{1}\right)  \right] \\
&  =\left[  \alpha_{2}\left(  C_{2}+D_{2}\right)  \right]  \otimes\left(
\alpha_{1}I\right)  +\left(  \alpha_{2}I\right)  \otimes\left[  \alpha
_{1}\left(  C_{1}+D_{1}\right)  \right] \\
&  =\mathbf{0}%
\end{align*}
and%
\[
\left(  \alpha_{2}\otimes\alpha_{1}\right)  \left(  \mathbf{e}\otimes
\mathbf{e}\right)  =\left(  \alpha_{2}\mathbf{e}\right)  \otimes\left(
\alpha_{1}\mathbf{e}\right)  =1.
\]
Therefore, $\alpha_{2}\otimes\alpha_{1}$ is the stationary probability vector
of the Markov process $\mathbb{A}_{k}$ for each $k\geq1$.

Now, we compute the (upward and downward) mean drift rates of the QBD process
$Q_{A}$. From Level $k$ to Level $k+1$, we obtain%
\begin{align*}
\left(  \alpha_{2}\otimes\alpha_{1}\right)  A_{0}^{\left(  k\right)
}\mathbf{e}  &  =\left(  \alpha_{2}\otimes\alpha_{1}\right)  \left(  I\otimes
D_{1}\right)  \left(  \mathbf{e}\otimes\mathbf{e}\right) \\
&  =\left(  \alpha_{2}I\mathbf{e}\right)  \otimes\left(  \alpha_{1}%
D_{1}\mathbf{e}\right)  =\lambda_{1}.
\end{align*}
Similarly, from Level $k$ to Level $k-1$, we get%
\begin{align*}
\text{\ }\left(  \alpha_{2}\otimes\alpha_{1}\right)  A_{2}^{\left(  k\right)
}\mathbf{e}  &  =\left(  \alpha_{2}\otimes\alpha_{1}\right)  \left[
D_{2}\otimes I+k\theta_{1}I\right]  \mathbf{e}\\
&  =\left(  \alpha_{2}\otimes\alpha_{1}\right)  \left(  D_{2}\otimes I\right)
\mathbf{e+}\left(  \alpha_{2}\otimes\alpha_{1}\right)  k\theta_{1}%
I\mathbf{e}\\
&  =\lambda_{2}+k\theta_{1}.
\end{align*}
Since $k$ is a positive integer and $\theta_{1}>0$, it is easy to check that
if $k>\max\left\{  1,\left(  \lambda_{1}-\lambda_{2}\right)  /\theta
_{1}\right\}  $, then $\left(  \alpha_{2}\otimes\alpha_{1}\right)
A_{0}^{\left(  k\right)  }\mathbf{e<}$ $\left(  \alpha_{2}\otimes\alpha
_{1}\right)  A_{2}^{\left(  k\right)  }\mathbf{e}$.\textbf{ }Therefore, the
QBD process $Q_{A}$ must be positive recurrent due to the fact that the
positive integer $k$ goes to infinity.

On the other hand, we can similarly discuss the stability of the QBD process
$Q_{B}$. Let $\mathbb{B}_{k}=B_{0}^{\left(  k\right)  }+B_{1}^{\left(
k\right)  }+B_{2}^{\left(  k\right)  }$. Then we obtain%
\begin{align*}
\mathbb{B}_{k}  &  =D_{2}\otimes I+C_{2}\oplus C_{1}+k\theta_{2}I+I\otimes
D_{1}-k\theta_{2}I\\
&  =\left(  C_{2}+D_{2}\right)  \oplus\left(  C_{1}+D_{1}\right)  ,
\end{align*}
which is independent of the negative integer $k\leq-1$, and $\mathbb{B}_{k}$
is the infinitesimal generator of a continuous-time Markov process with
$m_{1}m_{2}$ states.

For $k\leq-1$, we obtain%
\begin{align*}
\left(  \alpha_{2}\otimes\alpha_{1}\right)  \mathbb{B}_{k}  &  =\left(
\alpha_{2}\otimes\alpha_{1}\right)  \left[  \left(  C_{2}+D_{2}\right)
\oplus\left(  C_{1}+D_{1}\right)  \right] \\
&  =\left(  \alpha_{2}\otimes\alpha_{1}\right)  \left[  \left(  C_{2}%
+D_{2}\right)  \otimes I\right]  +\left(  \alpha_{2}\otimes\alpha_{1}\right)
\left[  I\otimes\left(  C_{1}+D_{1}\right)  \right] \\
&  =\left[  \alpha_{2}\left(  C_{2}+D_{2}\right)  \right]  \otimes\left(
\alpha_{1}I\right)  +\left(  \alpha_{2}I\right)  \otimes\left[  \alpha
_{1}\left(  C_{1}+D_{1}\right)  \right] \\
&  =\mathbf{0}%
\end{align*}
and%
\[
\left(  \alpha_{2}\otimes\alpha_{1}\right)  \left(  \mathbf{e}\otimes
\mathbf{e}\right)  =\left(  \alpha_{2}\mathbf{e}\right)  \otimes\left(
\alpha_{1}\mathbf{e}\right)  =1.
\]
Therefore, $\alpha_{2}\otimes\alpha_{1}$ is the stationary probability vector
of the Markov process $\mathbb{B}_{k}$ for each $k\leq-1$.

Now, we compute the (upward and downward) mean drift rates of the QBD process
$Q_{B}$. From Level $k$ to Level $k+1$, we yield%
\begin{align*}
\text{\ }\left(  \alpha_{2}\otimes\alpha_{1}\right)  B_{2}^{\left(  k\right)
}\mathbf{e}  &  =\left(  \alpha_{2}\otimes\alpha_{1}\right)  \left[  \left(
I\otimes D_{1}\right)  -k\theta_{2}I\right]  \mathbf{e}\\
&  =\left(  \alpha_{2}I\mathbf{e}\right)  \otimes\left(  \alpha_{1}%
D_{1}\mathbf{e}\right)  -\left(  \alpha_{2}\otimes\alpha_{1}\right)
k\theta_{2}I\mathbf{e}\\
&  =\lambda_{1}+\left(  -k\right)  \theta_{2}.
\end{align*}
Similarly, from Level $k$ to Level $k-1$, we have%
\begin{align*}
\left(  \alpha_{2}\otimes\alpha_{1}\right)  B_{0}^{\left(  k\right)
}\mathbf{e}  &  =\left(  \alpha_{2}\otimes\alpha_{1}\right)  \left(
D_{2}\otimes I\right)  \left(  \mathbf{e}\otimes\mathbf{e}\right) \\
&  =\left(  \alpha_{2}D_{2}\mathbf{e}\right)  \otimes\left(  \alpha
_{1}I\mathbf{e}\right)  =\lambda_{2}.
\end{align*}
Since $k$ is a negative integer and $\theta_{2}>0$, it is easy to check that
if $k<\max\left\{  -1,-\left(  \lambda_{2}-\lambda_{1}\right)  /\theta
_{2}\right\}  $, then $\left(  \alpha_{2}\otimes\alpha_{1}\right)
B_{0}^{\left(  k\right)  }\mathbf{e<}\left(  \alpha_{2}\otimes\alpha
_{1}\right)  B_{2}^{\left(  k\right)  }\mathbf{e}$. Therefore, the QBD process
$Q_{B}$ must be positive recurrent due to the fact that the negative integer
$k$ goes to infinity.

Based on the above analysis, the two QBD processes $Q_{A}$ and $Q_{B}$ are all
positive recurrent, so that the bilateral QBD process $Q$ is irreducible and
positive recurrent. Therefore, the block-structured double-ended queue is
stable. This completes the proof. \hfill$\blacksquare$

\vskip                                  0.4cm

\textbf{(b) Proof of Theorem \ref{The2}}

The proof is easy through checking that $\pi$ satisfies the system of linear
equations $\pi Q=\mathbf{0}$ and $\pi\mathbf{e}=1$.\ To this end, we consider
the following three different cases:

\textit{Case one:} $k\geq2$. In this case, we need to check that%
\begin{equation}
\pi_{k-1}A_{0}^{\left(  k-1\right)  }+\pi_{k}A_{1}^{\left(  k\right)  }%
+\pi_{k+1}A_{2}^{\left(  k+1\right)  }=0. \label{equ1211}%
\end{equation}
Since $\pi_{k-1}=c\widetilde{\pi}_{k-1}=c\widetilde{\pi}_{1}R_{1}R_{2}\cdots
R_{k-2},\pi_{k}=c\widetilde{\pi}_{k}=c\widetilde{\pi}_{1}R_{1}R_{2}\cdots
R_{k-1}$ and $\pi_{k+1}=c\widetilde{\pi}_{k+1}=c\widetilde{\pi}_{1}R_{1}%
R_{2}\cdots R_{k}$, we obtain%
\begin{align*}
&  \text{ \ \ }\pi_{k-1}A_{0}^{\left(  k-1\right)  }+\pi_{k}A_{1}^{\left(
k\right)  }+\pi_{k+1}A_{2}^{\left(  k+1\right)  }\\
&  =c\widetilde{\pi}_{1}R_{1}R_{2}\cdots R_{k-2}A_{0}^{\left(  k-1\right)
}+c\widetilde{\pi}_{1}R_{1}R_{2}\cdots R_{k-1}A_{1}^{\left(  k\right)
}+c\widetilde{\pi}_{1}R_{1}R_{2}\cdots R_{k}A_{2}^{\left(  k+1\right)  }\\
&  =c\widetilde{\pi}_{1}R_{1}R_{2}\cdots R_{k-2}\left(  A_{0}^{\left(
k-1\right)  }+R_{k-1}A_{1}^{\left(  k\right)  }+R_{k-1}R_{k}A_{2}^{\left(
k+1\right)  }\right)  =0
\end{align*}
by means of (\ref{equ3}).

\textit{Case two:} $k\leq-2$. In this case, we need to check that
\begin{equation}
\pi_{k+1}B_{0}^{\left(  k+1\right)  }+\pi_{k}B_{1}^{\left(  k\right)  }%
+\pi_{k-1}B_{2}^{\left(  k-1\right)  }=0. \label{equ1212}%
\end{equation}
Note that $\pi_{k+1}=c\widetilde{\pi}_{k+1}=c\widetilde{\pi}_{-1}%
\mathbb{R}_{-1}\mathbb{R}_{-2}\cdots\mathbb{R}_{k+2},\pi_{k}=c\widetilde{\pi
}_{k}=c\widetilde{\pi}_{-1}\mathbb{R}_{-1}\mathbb{R}_{-2}\cdots\mathbb{R}%
_{k+1}$ and $\pi_{k-1}=c\widetilde{\pi}_{k-1}=c\widetilde{\pi}_{-1}%
\mathbb{R}_{-1}\mathbb{R}_{-2}\cdots\mathbb{R}_{k}$, we have%
\begin{align*}
&  \text{ \ \ }\pi_{k+1}B_{0}^{\left(  k+1\right)  }+\pi_{k}B_{1}^{\left(
k\right)  }+\pi_{k-1}B_{2}^{\left(  k-1\right)  }\\
&  =c\widetilde{\pi}_{-1}\mathbb{R}_{-1}\mathbb{R}_{-2}\cdots\mathbb{R}%
_{k+2}B_{0}^{\left(  k+1\right)  }+c\widetilde{\pi}_{-1}\mathbb{R}%
_{-1}\mathbb{R}_{-2}\cdots\mathbb{R}_{k+1}B_{1}^{\left(  k\right)
}+c\widetilde{\pi}_{-1}\mathbb{R}_{-1}\mathbb{R}_{-2}\cdots\mathbb{R}_{k}%
B_{2}^{\left(  k-1\right)  }\\
&  =c\widetilde{\pi}_{-1}\mathbb{R}_{-1}\mathbb{R}_{-2}\cdots\mathbb{R}%
_{k+2}\left(  B_{0}^{\left(  k+1\right)  }+\mathbb{R}_{k+1}B_{1}^{\left(
k\right)  }+\mathbb{R}_{k+1}\mathbb{R}_{k}B_{2}^{\left(  k-1\right)  }\right)
=0
\end{align*}
in terms of (\ref{equ7}).

\textit{Case three:} $k=1,$ $0,$ $-1.$ In this case, we can check that%
\begin{equation}
\left\{
\begin{array}
[c]{l}%
\pi_{0}A_{0}^{\left(  0\right)  }+\pi_{1}A_{1}^{\left(  1\right)  }+\pi
_{2}A_{2}^{\left(  2\right)  }=0,\\
\pi_{-1}B_{2}^{\left(  -1\right)  }+\pi_{0}\left(  B_{1}^{\left(  0\right)
}+A_{1}^{\left(  0\right)  }\right)  +\pi_{1}A_{2}^{\left(  1\right)  }=0,\\
\pi_{-2}B_{2}^{\left(  -2\right)  }+\pi_{-1}B_{1}^{\left(  -1\right)  }%
+\pi_{0}B_{0}^{\left(  0\right)  }=0.
\end{array}
\right.  \label{equ1213}%
\end{equation}
by means of $\pi_{0}=c\widetilde{\pi}_{0},\pi_{1}=c\widetilde{\pi}_{1},\pi
_{2}=c\widetilde{\pi}_{2}=c\widetilde{\pi}_{1}R_{1},\pi_{-1}=c\widetilde{\pi
}_{-1}$ and $\pi_{-2}=c\widetilde{\pi}_{-2}=c\widetilde{\pi}_{-1}%
\mathbb{R}_{-1}$.

Let $\mathbf{Z}$ be the set of all integers, i.e., $\mathbf{Z}=\left\{
\ldots,-2,-1,0,1,2,\ldots\right\}  $. Note that $\pi\mathbf{e}=1$, we compute%
\begin{equation}
\sum_{k\in\mathbf{Z}}\pi_{k}\mathbf{e}=1 \label{equ1214}%
\end{equation}
by means of $\pi_{k}=c\widetilde{\pi}_{k}=c\widetilde{\pi}_{1}R_{1}R_{2}\cdots
R_{k-1}$ for $k\geq2,\pi_{1}=c\widetilde{\pi}_{1},\pi_{0}=c\widetilde{\pi}%
_{0},\pi_{-1}=c\widetilde{\pi}_{-1}$ and $\pi_{k}=c\widetilde{\pi}%
_{k}=c\widetilde{\pi}_{-1}\mathbb{R}_{-1}\mathbb{R}_{-2}\cdots\mathbb{R}%
_{k+1}$ for $k\leq-2$. Thus we have%
\begin{align*}
1  &  =\sum_{k\in\mathbf{Z}}\pi_{k}\mathbf{e}=\sum_{k\leq-2}\pi_{k}%
\mathbf{e+}\pi_{-1}\mathbf{e+}\pi_{0}\mathbf{e+}\pi_{1}\mathbf{e+}\sum
_{k=2}^{\infty}\pi_{k}\mathbf{e}\\
\text{\ }  &  =\sum_{k\leq-2}c\widetilde{\pi}_{-1}\mathbb{R}_{-1}%
\mathbb{R}_{-2}\cdots\mathbb{R}_{k+1}\mathbf{e+}c\widetilde{\pi}%
_{-1}\mathbf{e+}c\widetilde{\pi}_{0}\mathbf{e+}c\widetilde{\pi}_{1}%
\mathbf{e+}\sum\limits_{k=2}^{\infty}c\widetilde{\pi}_{1}R_{1}R_{2}\cdots
R_{k-1}\mathbf{e}\\
&  =c\left(  \sum_{k\leq-2}\widetilde{\pi}_{-1}\mathbb{R}_{-1}\mathbb{R}%
_{-2}\cdots\mathbb{R}_{k+1}\mathbf{e}+\widetilde{\pi}_{-1}\mathbf{e}%
+\widetilde{\pi}_{0}\mathbf{e}+\widetilde{\pi}_{1}\mathbf{e}+\sum
\limits_{k=2}^{\infty}\widetilde{\pi}_{1}R_{1}R_{2}\cdots R_{k-1}%
\mathbf{e}\right)  ,
\end{align*}
which gives the positive constant $c$ in (\ref{equ121}). This completes the
proof.\hfill$\blacksquare$

\end{document}